\documentclass[a4paper,oneside,12pt]{article}
\pdfoutput=1
\usepackage{cite, graphicx, url, amsmath}
\usepackage{subfigure}

\title{\bf Stability of Fractional-Order Systems \\
with Rational Orders }

\author{Ivo~Petr\'a\v{s} \\[1ex]
        Institute of Control and Informatization of Production Processes \\
        BERG Faculty, Technical University of Ko\v{s}ice \\
        B. N\v{e}mcovej 3, 042 00 Ko\v{s}ice, Slovak Republic\\
        Tel./Fax: +421-55-602-5194; E-mail: ivo.petras@tuke.sk\\}
\date{}

\begin{document}
\maketitle

\begin{abstract}
\noindent
This paper deals with stability of a certain class of fractional order
linear and nonlinear systems. The stability is investigated in the time domain and the frequency domain.
The general stability conditions and several illustrative examples are presented as well.
\end{abstract}


\section{Introduction}

Fractional calculus is more than 300 years old topic.
A number of applications where fractional calculus has been used rapidly
grows. This mathematical phenomena allow to describe a real object more
accurate that the classical ``integer'' methods. The real objects are
generally fractional \cite{Oustaloup, Podlubny, Nakagava, Westerlund2, Wang},
however, for many of them the fractionality is very low.
The main reason for using the integer-order models was the absence of
solution methods for fractional differential equations.

Recently, the fractional order linear time invariant (FOLTI)
systems have attracted lots of attention in control
systems society (e.g.: \cite{Lurie, Dorcak, Podlubny-PID, Oustaloup, Raynaud}) even
though fractional-order control problems were investigated
as early as 1960's \cite{Manabe}. In the fractional order controller,
the fractional order integration or derivative of the output
error is used for the current control force calculation.

The fractional order calculus plays an important role in physics \cite{Parada, Tapia, Torvik},
thermodynamics \cite{Jesus, Vinagre4}, electrical circuits theory and fractances \cite{Arena2, Bode, Charef, Nakagava, Carlson, Westerlund},
mechatronics systems \cite{Silva}, signal processing \cite{Vinagre, Tseng}, chemical mixing \cite{Oldham}, chaos theory \cite{Haeri, Tavazoei4},
and biological system as well \cite{Chen06}. It is recommended to refer to (e.g.: \cite{Oustaloup2,
Axtell, Vinagre6, Xue, Monje-cep08}) for the further engineering applications
of fractional order systems. The question of stability is very important especially in control theory.
In the field of fractional-order control systems,
there are many challenging and unsolved problems related to stability theory such as
robust stability, bounded input - bounded output stability, internal stability, root-locus,
robust controllability, robust observability, etc.

For distributed parameter systems with a distributed delay \cite{Ozturk}, provided an stability analysis
method which may be used to test the stability of fractional order
differential equations.  In \cite{Bonnet00}, the co-prime
factorization method is used for stability analysis of fractional
differential systems. In \cite{Matignon96}, the stability conditions
for commensurate FOLTI system have been provided.
However, the general robust stability test procedure
and proof of the validity for the general type of the FOLTI system is still
open and discussed in \cite{Petras4}.
Stability has also been investigated for fractional order nonlinear system
(chaotic system) with commensurate and incomensurate order as well
\cite{Ahmed, Haeri, Tavazoei2}.

This paper is organized as follows. In Sec.~2 is briefly introduced the fractional calculus.
Sec.~3 is on fractional order systems. In Sec.~4 are analyzed the stability conditions of fractional
order linear and nonlinear systems. Sec.~5 concludes this paper with some remarks.

\section{Fractional Calculus Fundamentals}

\subsection{Definitions of Fractional Derivatives and Integrals}

The idea of fractional calculus has been known since the development
of the regular calculus, with the first reference
probably being associated with Leibniz and L'Hospital in 1695
where half-order derivative was mentioned.

Fractional calculus is a generalization of integration and differentiation
to non-integer order fundamental operator $_{a}D^{r}_{t}$,
where $a$ and $t$ are the limits of the operation and $r \in R$.
The continuous integro-differential operator is defined as
$$
_aD^{r}_{t} = \left \{
        \begin{array}{ll}
                \frac{d^{r}}{dt^{r}} & : r>0, \\
                 1 & : r=0, \\
                \int_{a}^{t} (d\tau)^{-r} & : r<0.
        \end{array}
        \right.
$$

The three definitions used for the general fractional differintegral
are the Grunwald-Letnikov (GL) definition,  the Riemann-Liouville (RL)
and the Caputo definition \cite{Oldham, Podlubny}. The GL is given here
\begin{equation}\label{GLD_d}
     _{a}D^{r}_{t}f(t)=\lim_{h \to 0} h^{-r}
     \sum_{j=0}^{[\frac{t-a}{h}]}(-1)^j {r \choose j}f(t-jh),
\end{equation}
where $[.]$ means the integer part. The RL definition is given as
\begin{equation}\label{LRL}
   _{a}D_{t}^{r}f(t)=
    \frac{1}{\Gamma (n - r)}
    \frac{d^{n}}{dt^{n}}
    \int_{a}^{t}
    \frac{f(\tau)}{(t-\tau)^{r - n + 1}}d\tau,
\end{equation}
for $(n-1 < r <n)$ and  where $\Gamma (.)$ is the {\it Gamma}  function.
The Caputo's definition can be written as
\begin{equation}\label{Caputo}
   _{a}D_{t}^{r}f(t)=
    \frac{1}{\Gamma (r - n)}
    \int_{a}^{t}
    \frac{f^{(n)}(\tau)}{(t-\tau)^{r - n + 1}}d\tau,
\end{equation}
for $(n-1 < r <n)$. The initial conditions for the fractional order differential equations with the Caputo's derivatives
are in the same form as for the integer-order differential equations.

\newpage
\subsection{Some Properties of Fractional Derivatives and Integrals}

The main properties of fractional derivatives and integrals
are the following:

\begin{enumerate}

\item
If $f(t)$ is an analytical function of $t$, then its  fractional
derivative $_{0}D_{t}^{\alpha}f(t)$ is an analytical function of
$t$, $\alpha$.

\item
For $\alpha=n$, where $n$ is integer, the operation $_{0}D_{t}^{\alpha}f(t)$
gives the same result as classical differentiation of integer order $n$.

\item
For $\alpha=0$ the operation $_{0}D_{t}^{\alpha}f(t)$ is the identity
operator:
$$_{0}D_{t}^{0}f(t) = f(t)$$

\item
Fractional differentiation and fractional integration are linear operations:
$$
_aD^r_t \left( \lambda f(t) + \mu g (t) \right )= \lambda\,
_aD^r_t f(t) + \mu\, _aD^r_t g(t).
$$

\item
The additive index law (semigroup property)
$$
_{0}D_{t}^{\alpha} \, _{0}D_{t}^{\beta} f(t) = \,
_{0}D_{t}^{\beta} \, _{0}D_{t}^{\alpha} f(t) = \,
 _{0}D_{t}^{\alpha + \beta} f(t)
$$
holds under some reasonable constraints on the function $f(t)$.

The fractional-order derivative commutes with integer-order derivation
$$
\frac{d^n}{dt^n} (_aD^r_t f(t))= \,_aD^r_t \left ( \frac{d^n f(t)}{dt^n}\right ) =
	\, _aD_t^{r+n}f(t),
$$
under the condition $t=a$ we have
$f^{(k)}(a)=0,\,(k=0, 1, 2,\dots,n-1)$.
The relationship above says the operators $\frac{d^n}{dt^n}$ and $_aD_t^r$ commute.

\item
The formula for the Laplace transform of the RL fractional derivative
(\ref{LRL}) has the form \cite{Podlubny}:
$$
      \int_{0}^{\infty} e^{-st}\, _{0}D_{t}^{r}f(t) \, dt =
      s^{r}F(s) - \sum_{k=0}^{n-1} s^{k} \,
                    \left.  _{0}D_{t}^{r -k-1}f(t) \right|_{t=0},
$$
for $(n-1 < r \leq n)$, where $s \equiv j\omega$ denotes the Laplace operator.
For zero initial conditions, Laplace transform of fractional derivatives
(Grunwald-Letnikov, Riemann-Liouville, and Caputo's), reduces to:
$$
  L \{_0D^r_t f(t)\} = s^r F(s).
$$

\item   Geometric and physical interpretation of fractional integration
and fractional differentiation were exactly described in Podlubny's work
\cite{Podlubny-geomet}.

\end{enumerate}

Some others important properties of the fractional derivatives and integrals
as for example Leibniz's rule, translation, Chain rule, bahaviour and dependence on limit and so on, we can find out
in several works (e.g.: \cite{Oldham, Podlubny, Oustaloup}, etc.).

\section{Fractional-Order Systems}

\subsection{Fractional LTI Systems}

A general fractional-order system can be described by a fractional differential equation of the form
\begin{eqnarray}
a_{n}D^{\alpha _{n}}y(t)+a_{n-1}D^{\alpha _{n-1}}y(t)+\ldots +a_{0}D^{\alpha
_{0}}y(t) = && \hspace*{4em} \nonumber\\
&& \hspace*{-18em}
 = b_{m}D^{\beta _{m}}u(t)+b_{m-1}D^{\beta _{m-1}}u(t)+\ldots +b_{0}D^{\beta
_{0}}u(t),
\end{eqnarray}
or by the corresponding transfer function of \textit{incommensurate} real orders of the following
form \cite{Podlubny}:
\begin{equation}
G(s)=\frac{b_{m}s^{\beta _{m}}+\ldots +b_{1}s^{\beta
_{1}}+b_{0}s^{\beta _{0}}}{%
a_{n}s^{\alpha _{n}}+\ldots +a_{1}s^{\alpha _{1}}+a_{0}s^{\alpha
_{0}}}=\frac{%
Q(s^{\beta _{k}})}{P(s^{\alpha _{k}})},  \label{Eq1}
\end{equation}
where $D^{\gamma} \equiv \,_{0}D_{t}^{\gamma}$ denotes the Riemann-Liouville
or Caputo fractional derivative \cite{Podlubny};
$a_k$ $(k=0,\ldots\,n)$, $b_k$ $(k=0,\ldots\,m)$ are
constant; and
$\alpha_k$ $(k=0,\ldots\,n)$, $\beta_k$ $(k=0,\ldots\,m)$
are arbitrary real numbers and without loss of generality they can be arranged as
$\alpha_n > \alpha_{n-1} > \ldots > \alpha_0$,
and
$\beta_m > \beta_{m-1} > \ldots > \beta_0$.

The incommensurate order system (\ref{Eq1}) can also be expressed in commensurate
form by the multi-valued transfer function \cite{Bayat2}
\begin{equation}\label{ITF}
   H(s) = \frac{b_m s^{m/v}+\dots +b_1 s^{1/v} +b_0}
               {a_n s^{n/v}+\dots +a_1 s^{1/v} +a_0}, \,\,\ (v>1).
\end{equation}
Note that every fractional order system can be expressed in the form (\ref{ITF}) and
domain of the $H(s)$ definition is a Riemann surface with $v$ Riemann sheets \cite{Lepage}.

In the particular case of \textit{commensurate} order systems, it
holds that,
$\alpha_{k}=\alpha k,\beta_{k}=\alpha k, (0<\alpha<1), \forall
k\in \mbox{Z}$, and the transfer function has the following form:
\begin{equation}
G(s)=K_{0} \frac{ \sum_{k=0}^{M} b_{k} (s^{\alpha})^{k} }
                { \sum_{k=0}^{N} a_{k} (s^{\alpha})^{k} }
    =K_{0} \frac{ Q (s^{\alpha})} { P(s^{\alpha})}  \label{Eq2}
\end{equation}
With $N>M$, the function $G(s)$ becomes a proper rational function
in the complex variable $s^{\alpha}$  which can be expanded in partial fractions of
the following form:
\begin{equation}
G(s)=K_{0} \left[ \sum_{i=1}^N  \frac{A_{i}}{s^{\alpha} +\lambda _{i}
 } \right],
\label{Eq3}
\end{equation}
where  $\lambda _{i}\, (i=1,2,..,N)$ are the roots of the pseudo-polynomial
$P(s^{\alpha})$ or the system poles which are assumed to be simple without loss of
generality.  The analytical solution of the system (\ref{Eq3}) can
be expressed as
\begin{equation}
   y(t) = L^{-1} \left \{ K_{0} \left[ \sum_{i=1}^N  \frac{A_{i}}{s^{\alpha} +\lambda _{i}} \right] \right \}
   = K_0 \sum_{i=1}^{N} A_i t^{\alpha} E_{\alpha,\alpha}(-\lambda_i t^{\alpha}).
\end{equation}

A fractional order plant to be controlled can be described by a typical $n$-term linear
homogeneous fractional order differential equation (FODE) in time domain
\begin{equation}
a_{n}\, D^{\alpha_{n}}_t y(t) + \cdots + a_{1}\,D
^{\alpha_{1}}_t y(t) + a_{0}\, D^{\alpha_{0}}_t y(t) = 0
\label{n-term-equation}
\end{equation}
where $a_k (k= 0, 1, \cdots, n)$ are constant coefficients of the FODE; $%
\alpha_k, (k = 0, 1, 2, \cdots, n)$ are real numbers. Without loss of
generality, assume that $\alpha_n > \alpha_{n-1} > \ldots > \alpha_0 \geq 0$.

The analytical solution of the FODE (\ref{n-term-equation}) is given by general formula \cite{Podlubny}
\begin{eqnarray} \label{Step-response}
     y (t)
     &=& \frac{1}{a_n} \sum_{m=0}^{\infty}
                     \frac{(-1)^m}{m!}
			   \hspace*{-3mm}
                     \sum_{{k_0+k_1+\ldots +k_{n-2} = m  \atop
                     k_0 \geq 0; \ldots\, , k_{n-2} \geq 0  }}
			   \hspace*{-7mm}
                     (m; k_0, k_1, \ldots \,, k_{n-2}) \nonumber \\
     &\times& \prod_{i=0}^{n-2}
     \left(
          \frac{a_i}{a_n}
     \right)^{k_i}
     {\cal E}_{m}(t, -\frac{a_{n-1}}{a_{n}}; \alpha_{n}-\alpha_{n-1},
                      \alpha_n \nonumber \\
     &+& \sum_{j=0}^{n-2}(\alpha_{n-1}-\alpha_j)k_j + 1),
     \hspace*{-5mm}
\end{eqnarray}
where $(m; k_0, k_1, \ldots \,, k_{n-2})$ are the multinomial coefficients
and ${\cal E}_{k}(t, y; \mu, \nu)$ is the function of Mittag-Leffler type
introduced by Podlubny \cite{Podlubny}. The function is defined by
\begin{equation} \label{Podlubny-function}
      {\cal E}_{k}(t, y; \mu, \nu) =
      t^{\mu k +\nu -1} E_{\mu, \nu}^{(k)}(yt^{\mu}),
      \quad
      (k = 0, 1, 2, \ldots),
\end{equation}
where $E_{\mu, \nu}(z)$ is the Mittag-Leffler function of two parameters \cite{Gorenflo}:
\begin{equation} \label{ML-Definition-Control}
E_{\mu, \nu}(z) = \sum_{i=0}^{\infty}\frac{z^i}{\Gamma (\mu i + \nu)},
\hspace{3em} (\mu > 0, \hspace{1em} \nu > 0),
\end{equation}
where e.g. $E_{1,1}(z) = e^{z}$,
and where its $k$-th derivative is given by
\begin{equation}
E_{\mu,\nu}^{(k)}(z) =
\sum_{i=0}^{\infty} \frac{(i+k)! \,\, z^{i}}
                         {i! \,\, \Gamma (\mu i + \mu k + \nu)},
\hspace{3em}
(k = 0, 1, 2, ...).
\end{equation}
Consider a control function which acts on the FODE system
(\ref{n-term-equation}) as follows:
\begin{equation}
a_{n}\, D^{\alpha_{n}}_ty(t) + \cdots + a_{1}\, D^{\alpha_{1}}_ty(t) + a_{0}\, D^{\alpha_{0}}_ty(t) = u(t).
\label{n-term-equation_u}
\end{equation}
By Laplace transform, we can get a fractional transfer function:
\begin{equation}  \label{Gplant}
G(s) =\frac{Y(s)}{U(s)}= \frac{1}{a_{n}s^{\alpha_{n}} + \cdots +
a_{1}s^{\alpha_{1}} + a_{0}s^{\alpha_{0}}}.
\end{equation}
The fractional order linear time-invariant system can also be represented by the following state-space model
\begin{eqnarray}\label{LTI_SS}
_0D^{\textbf{q}}_{t} x(t) &=& \textbf{A} x(t) + \textbf{B} u(t) \nonumber \\
y(t) &=& \textbf{C} x(t)
\end{eqnarray}
where $x\in R^n$, $u \in R^r$ and $y \in R^p$ are the state, input and output vectors of the system and $\textbf{A} \in R^{n\times n}$,
$\textbf{B} \in R^{n\times r}$, $\textbf{C}\in R^{p\times n}$, and $\textbf{q}$ = $[q_1, q_2, \dots, q_n]^T$ are the fractional orders.
If $q_1=q_2=\dots q_n$, system (\ref{LTI_SS}) is called a commensurate order system, otherwise it is an incommensurate order system.

A fractional-order system described by $n$-term fractional differential equation (\ref{n-term-equation_u})
can be rewritten to the state-space representation in the form \cite{Dorcak1, Yang}:
\begin{eqnarray} \label{LTI_SS_M}
\left[
\begin{array}{c}
_0D^{q_1 }x_{1}(t) \\
_0D^{q_2 }x_{2}(t) \\
. \\
. \\
_0D^{q_n }x_{n}(t)
\end{array}
\right] &=&\left[
\begin{array}{ccccc}
0 & 1 & . & . & 0 \\
0 & 0 & 1 & . & 0 \\
. & . & . & . & . \\
. & . & . & . & . \\
-a_0/a_n & -a_1/a_n & . & . & a_{n-1}/a_n
\end{array}
\right] \left[
\begin{array}{c}
x_{1}(t) \\
x_{2}(t) \\
. \\
. \\
x_{n}(t)
\end{array}
\right] +\left[
\begin{array}{c}
0 \\
0 \\
. \\
. \\
1/a_n
\end{array}
\right] u(t) \nonumber \\
y(t) &=&\left[
\begin{array}{cccccccc}
1 & 0 & . & . & . & . & 0 & 0
\end{array}
\right] \left[
\begin{array}{c}
x_{1}(t) \\
x_{2}(t) \\
. \\
. \\
x_{n}(t)
\end{array}
\right],
\end{eqnarray}
where $\alpha_0=0$, $q_1=\alpha_1$, $q_2=\alpha_{n-1}-\alpha_{n-2}$, \dots $q_n=\alpha_n-\alpha_{n-1}$, and
with initial conditions:
\begin{eqnarray}\label{InitCond}
x_1(0) &=& x_{0}^{(1)}=y_0, \,\,\, x_2(0)=x_{0}^{(2)}=0, \dots \nonumber \\
x_i(0) &=& x_{0}^{(i)}=
	\left \{
        \begin{array}{ll}
                y_0^{(k)}, & \mbox{if}\,\,\, i=2k+1, \\
                0,          & \mbox{if}\,\,\, i=2k,
        \end{array}
        \right. \,\, i \leq n.
\end{eqnarray}

The $n$-term FODE (\ref{n-term-equation_u}) is equivalent to the system of equations (\ref{LTI_SS_M})
with the initial conditions (\ref{InitCond}).

Similar to conventional observability and controllability concept, the controllability
is defined as follow \cite{MatignonNovel}: System (\ref{LTI_SS}) is \textit{controllable} on $[t_0, t_{final}]$
if controllability matrix $C_a = [B|AB|A^2B|\dots|A^{n-1}B]$ has rank $n$. The observability is
defined as follow \cite{MatignonNovel}: System (\ref{LTI_SS}) is \textit{observable} on $[t_0, t_{final}]$
if observanility matrix $O_a = [C|CA|CA^2|\dots|CA^{n-1}]^T$ has rank $n$.

\subsection{Fractional Nonlinear Systems}

Generally, we consider the following incommensurate fractional order nonlinear system in the form:
\begin{eqnarray} \label{FONS}
     _0D^{q_i}_t x_i (t) &=& f_i (x_1(t), x_2(t), \dots, x_n(t), t) \nonumber \\
                 x_i (0) &=& c_i, \,\,\, i=1, 2, \dots, n,
\end{eqnarray}
where $c_i$ are initial conditins, or in its vector representation:
\begin{equation}\label{general}
    D^{\textbf{q}} \textbf{x} = \textbf{f}(\textbf{x}),
\end{equation}
where $\textbf{q}=[q_1, q_2, \dots, q_n]^T$ for $0<q_i<2$, $(i=1, 2, \dots, n)$ and $\textbf{x} \in R^n$.

The equilibrium points of system (\ref{general}) are
calculated via solving the following equation
\begin{equation}
\textbf{f}(\textbf{x}) = 0
\end{equation}
and we suppose that $x^{*} = (x_1^{*}, x_2^{*}, \dots, x_n^{*})$ is an equilibrium point of system (\ref{general}).

\section{Stability of the Fractional Order Systems}

\subsection{Preliminary Consideration}

Stability as an extremely important property of the dynamical systems
can be investigated in various domain \cite{Dorf-Bishop, Dazzo}. Usual concept of bounded input - bounded output (BIBO) or external
stability in \textit{time domain} can be defined via the following general stability conditions \cite{Matignon}:

A causal LTI system with impulse response $h(t)$ to be BIBO stable if the necessary
and sufficient condition is satisfied
$$
\int_0^{\infty} ||h(\tau)|| d\tau < \infty,
$$
where output of the system is defined by convolution
$$
  y(t)=h(t)*u(t) = \int_0^{\infty} h(\tau) u(t-\tau) d\tau,
$$
where $u, y \in L_{\infty}$ and $h \in L_1$.

Another very important domain is \textit{frequency domain}. In the case of frequency method
for evaluating the stability we transform the $s$-plane into the complex plane
$G_o(j\omega)$ and the transformation is realized according to the transfer function of the open loop
system $G_o(j\omega)$. During the transformation, all roots of the characteristic polynomial are
mapped from $s$-plane into the critical point $(-1, j0)$ in the plane $G_o(j\omega)$.
The mapping of the $s$-plane into $G_o(j\omega)$ plane is conformal, that is, the direction
and location of points in the $s$-plane is preserved in the $G_o(j\omega)$ plane.
Frequency investigation method and utilization of the Nyquist frequency
characteristics based on argument principle were described in the paper \cite{Petras2}.

However, we can not directly use an algebraic tools as for example Routh-Hurwitz criteria
for the fractional order system because we do not have a characteristic polynomial
but pseudo-polynomial with rational power - \textit{multivalued function}. It is possible only
in some special cases \cite{Ahmed}. Moreover, modern control method as for example
LMI (Linear Matrix Inequality) methods \cite{Oustaloup2} or other algorithms \cite{Hamamci, Hwang2} already have been developed.
The advantage of LMI methods in control theory is due their connection with the
Lyapunov method (existence a quadratic Lyapunov function). More generally, LMI
methods are useful to test of matrix eigenvalues belong to a certain region in complex plane.
A simple test can be used \cite{Anderson}. Roots of polynomial $P(s)=\mbox{det}(s I - A)$ lie inside in region
$-\pi/2-\delta < \mbox{arg}  (s) < \pi/2 + \delta$ if eigenvalues of the matrix
\begin{eqnarray}  \label{matrixstability}
A_1 = \left[
\begin{array}{cc}
A\, \mbox{cos}\, \delta & -A\, \mbox{sin}\, \delta \\
A \, \mbox{sin}\, \delta & A \, \mbox{cos}\, \delta
\end{array}
\right]  \equiv
A \otimes \left[
\begin{array}{cc}
\mbox{cos}\, \delta & -\, \mbox{sin}\, \delta \\
\mbox{sin}\, \delta &  \, \mbox{cos}\, \delta
\end{array}
\right]
\end{eqnarray}
have negative real part, where $\otimes$ denotes Kronecker product. This property has been used to stability analysis of ordinary fractional
order LTI system and also for interval fractional order LTI system \cite{Tavazoei5}.

\begin{figure}[!ht]
\centering
\noindent
  \includegraphics[width=3in]{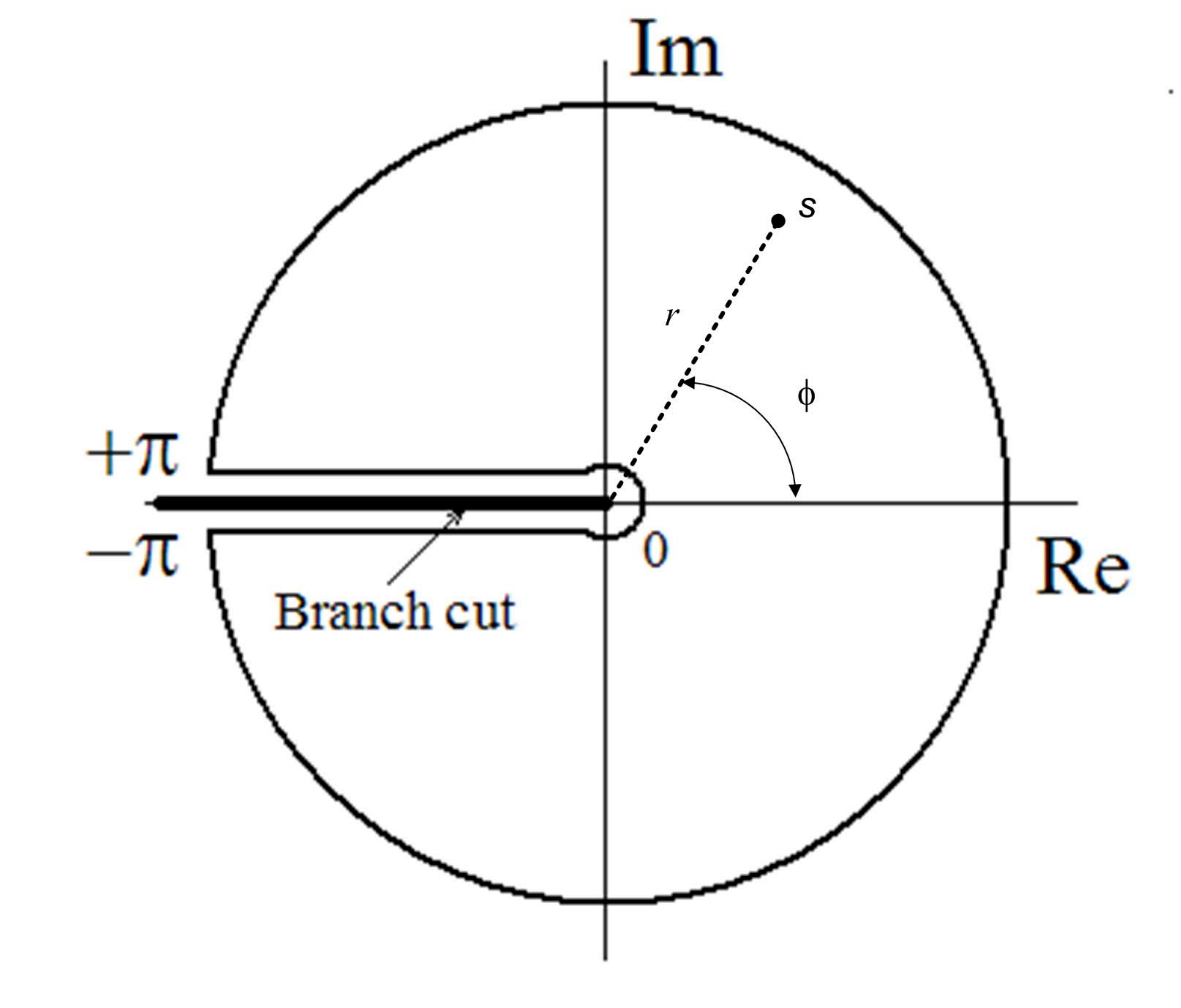}
   \caption{Branch cut $(0,-\infty)$ for branch points in the complex plane.}
   \label{Riemann1}
\end{figure}

When dealing with incommensurate fractional order systems (or, in general, with
fractional order systems) it is important to bear in mind that $P\left( s^{\alpha
}\right)$, $\alpha \in \mbox{R}$ is a multivalued function of $s^{\alpha}$,  $\alpha=\frac{u}{v}$,
the domain of which can be viewed as a Riemann surface with finite number
of Riemann sheets $v$, where origin is a branch point and the branch cut
is assumed at $\mbox{R}^-$ (see Fig.~\ref{Riemann1}).  Function  $s^{\alpha}$
becomes holomorphic in the complement of the branch cut line.
It is a fact that in multivalued functions only the
first Riemann sheet has its physical significance \cite{Gross}.
Note that each Riemann sheet has only one edge at branch cut and
not only poles and singularities originated from the characteristic
equation, but branch points and branch cut of given multivalued functions
are also important for the stability analysis \cite{Bayat4}.

In this paper the branch cut is assumed at $\mbox{R}^-$ and the first Riemann sheet is
denoted by $\Omega$ and defined as
\begin{equation}
          \Omega : =\{ r e^{j\phi}\, |\, r>0, -\pi < \phi < \pi \}.
\end{equation}

It is well-known that an integer order  LTI system is stable if all
the roots of the characteristic polynomial $P(s)$ are negative or have negative
real parts if they are complex conjugate (e.g.: \cite{Dorf-Bishop}). This
means that they are located on the left of the imaginary axis of the complex $s$-plane.
System $G(s)=Q(s)/P(s)$ is BIBO stable if
$$
\exists, \,\,\,\ ||G(s)|| \leq M < \infty, \,\,\,M>0, \,\,\,\, \forall s, \Re(s)\geq 0.
$$
A necessary and sufficient condition for the asymptotic stability is \cite{Salam}:
\begin{displaymath}
  \mbox{lim}_{t \to \infty} ||X(t)|| =0.
\end{displaymath}
According the final value theorem proposed in \cite{Ghartemani}, for fractional order case,
when there is a~branch point at $s=0$, we assume that $G(s)$ is multivalued function of
$s$, then
$$
  x(\infty) = \mbox{lim}_{s \to 0} [s G(s)].
$$

\begin{figure}[h]
\centering
\noindent
  \includegraphics[width=3.51in]{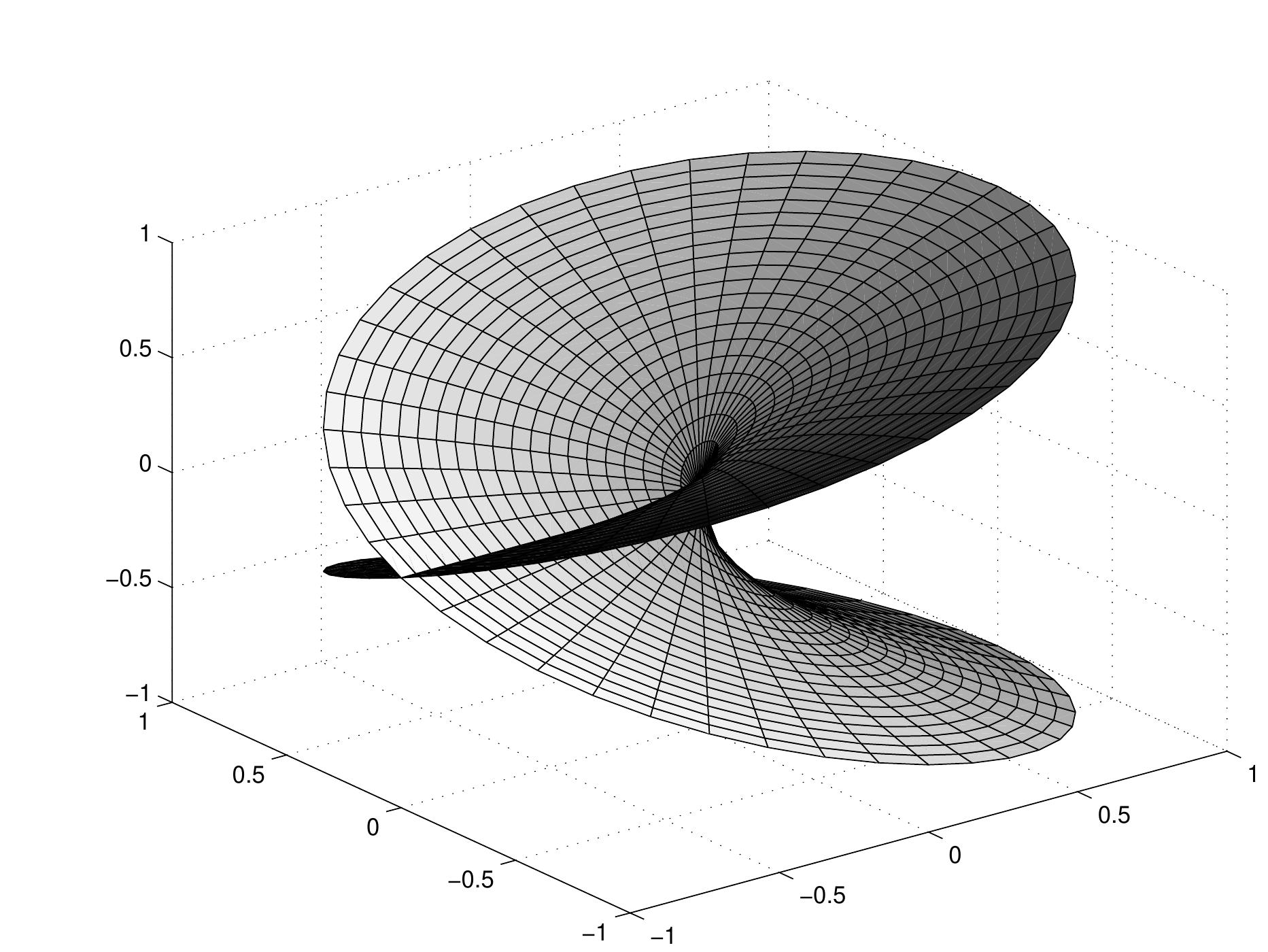}
   \caption{Riemann surface interpretation of the function $w=s^{\frac{1}{2}}$.}
   \label{RiemannS}
\end{figure}

\noindent
\textbf{Example 1}:
Let us investigate the simplest multi-valued function defined as follow
\begin{equation}\label{ss2}
w=s^{\frac{1}{2}}
\end{equation}
and there will be two $s$-planes which map onto a single $w$-plane. The interpretation
of the two sheets of the Riemann surface and the branch cut is depicted in Fig.~\ref{RiemannS}.

Define the principal square root function as
$$
f_1(s) = |s|^{\frac{1}{2}} e^{\frac{j\phi}{2}} = r e^{\frac{j\phi}{2}},
$$
where $r>0$ and $-\pi < \phi < +\pi$. The function $f_1(s)$ is a branch of $w$.
Using the same notation, we can find other branches of the square root function.  For example, if we let
$$
f_2(s) = |s|^{\frac{1}{2}} e^{\frac{j\phi + 2\pi}{2}} = r e^{\frac{j\phi + 2\pi}{2}},
$$
then $f_2(s) = - f_1(s)$ and it can be thought of as "plus" and "minus" square root functions.
The negative real axis is called a branch cut for the functions $f_1(s)$ and $f_2(s)$.
Each point on the branch cut is a point of discontinuity for both functions $f_1(s)$ and $f_2(s)$.
As has been shown in \cite{Lepage}, the function described by (\ref{ss2}) has a branch point of order 1 at $s=0$
and at infinity. They are located at ends of the branch cut (see also Fig.~\ref{Riemann1}).

\noindent
\textbf{Example 2}:
Let us investigate the transfer function of fractional-order system (multivalued function) defined as
\begin{equation}\label{sa}
G(s) = \frac{1}{s^{\alpha} + b},
\end{equation}
where $\alpha \in R$ $(0 < \alpha \leq 2)$ and $b \in R$ $(b>0)$.

The analytical solution of the fractional order system (\ref{sa}) obtained according to relation (\ref{Step-response})
has the following form:
\begin{equation}
g(t) =   {\cal E}_{0} (t, -b; \alpha, \alpha).
\end{equation}

The Riemann surface of the function (\ref{sa}) contains an infinite number of sheets
and infinitely many poles in positions
$$
  s = b^{\frac{1}{\alpha}} e^{\frac{j(\pi + 2\pi n)}{\alpha}}, \,\,\,\,\, n = 0, \pm 1, \pm 2, \ldots \,\,\,\,, \mbox{for}\,\,(\alpha>0)\,\mbox{and}\, (b>0).
$$
The sheets of the Riemann surface are all different if $\alpha$ is irrational.

For $1 < \alpha <2$ we have two poles corresponding to $n=0$ and $n=-1$, and poles are
$$
s = b^{\frac{1}{\alpha}} e^{\pm\frac{i\pi}{\alpha}}.
$$

\begin{figure}[!h]
\centering
\noindent
  \includegraphics[width=4.5in]{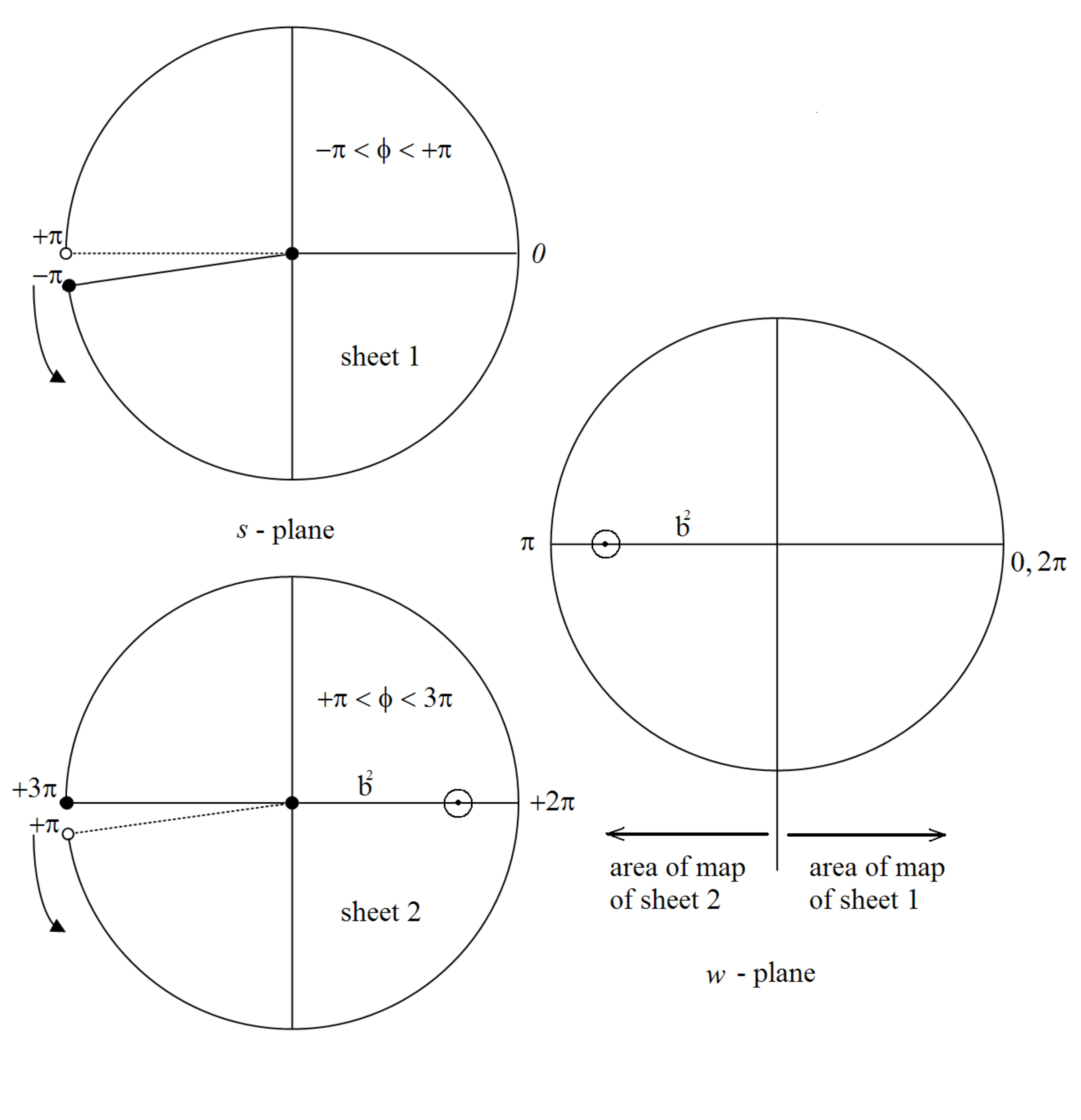}
   \caption{Correspondence between the $s$-plane and the $w$-plane for Eq.(\ref{sa22}).}
   \label{Riemann2}
\end{figure}

However, for $0< \alpha <1$ in (\ref{sa}) the denominator is a multivalued function and singularity of system can not
be defined unless it is made singlevalued. Therefore we will use the Riemann surface. Let us investigate transfer function
(\ref{sa}) for $\alpha=0.5$ (half-order system), then we get
\begin{equation}\label{sa22}
G(s) = \frac{1}{s^{\frac{1}{2}} + b},
\end{equation}
and by equating the denominator to zero we have
$$
s^{\frac{1}{2}} + b=0.
$$
Rewriting the complex operator $s^{\frac{1}{2}}$ in exponential form and using the well known relation
$e^{j\pi} +1 =0$ (or $e^{j(\pm \pi+2k \pi)} +1 =0$) we get the following formula:
\begin{equation}\label{sa2}
r^{\frac{1}{2}} e ^{j(\phi /2 +k \pi)} = a e^{j(\pm \pi + 2k \pi)}
\end{equation}
From relationship (\ref{sa2}) can be deduced that the modulus and phase (arg) of the pole are:
$$
r = b^2 \,\,\,\, \mbox{and} \,\,\, \phi=\pm2\pi(1+k)\,\,\, \mbox{for}\,\, k=0,1,2,\dots
$$
However the first sheet of the Riemann surface is defined for range of $-\pi < \phi < +\pi$,
the pole with the angle $\phi=\pm 2\pi$ does not fall within this range but pole
with the angle $\phi=2\pi$ falls to the range of the second sheet defined for $\pi < \phi < 2\pi$.
Therefore this half-order pole with magnitude $b^2$ is located on the second sheet of the Riemann
surface that consequently maps to the left side of the $w$-plane (see Fig.~\ref{Riemann2}). On this
plane the magnitude and phase of the singlevalued pole are $b^2$ and $\pi$, respectively \cite{Lepage}.

\begin{figure}[h]
\centering
\noindent
  \subfigure[Riemann surface]{
		\includegraphics[width=3.5in]{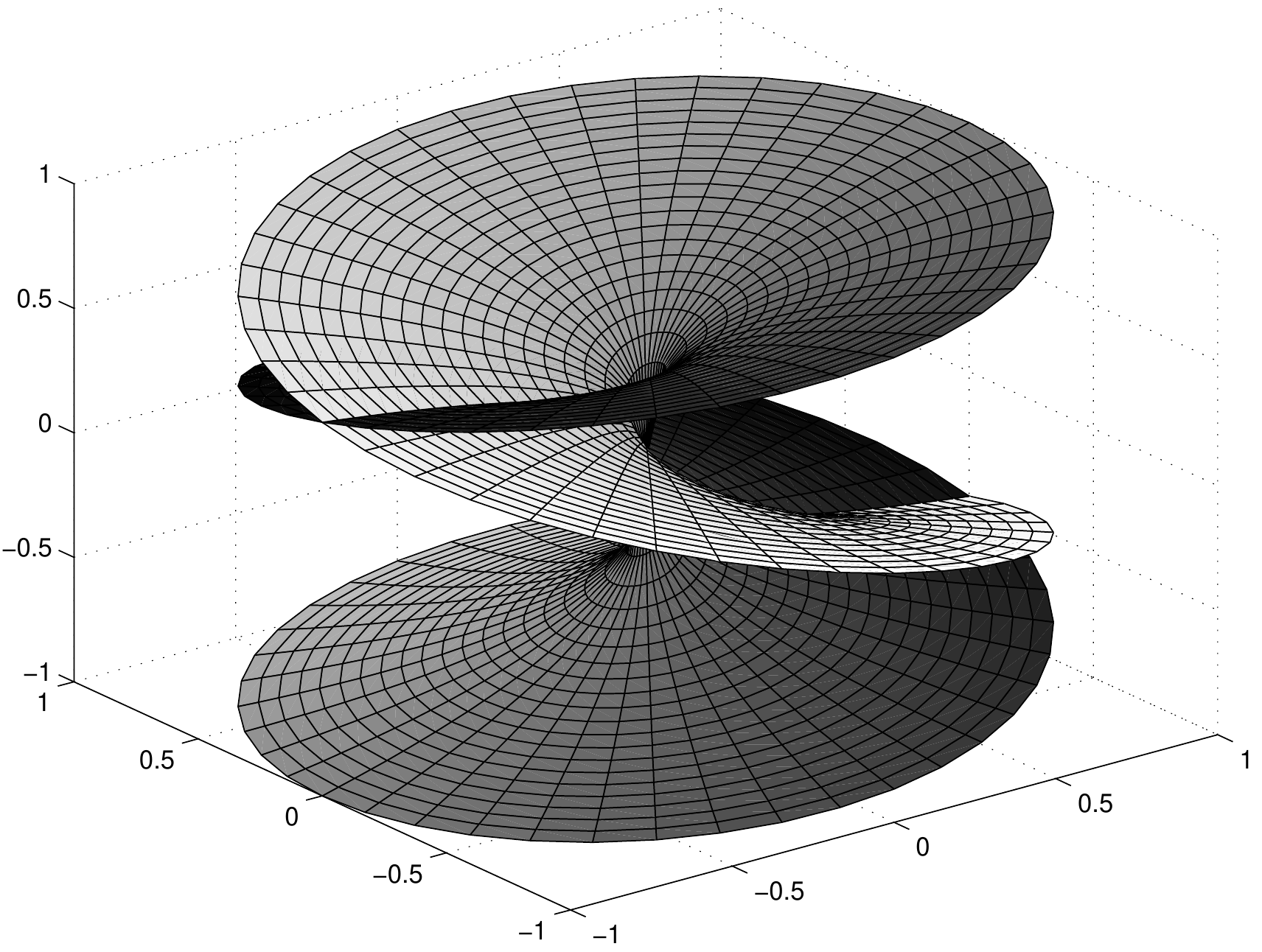}}
  \subfigure[Complex $w$-plane]{
                \includegraphics[width=2.9in]{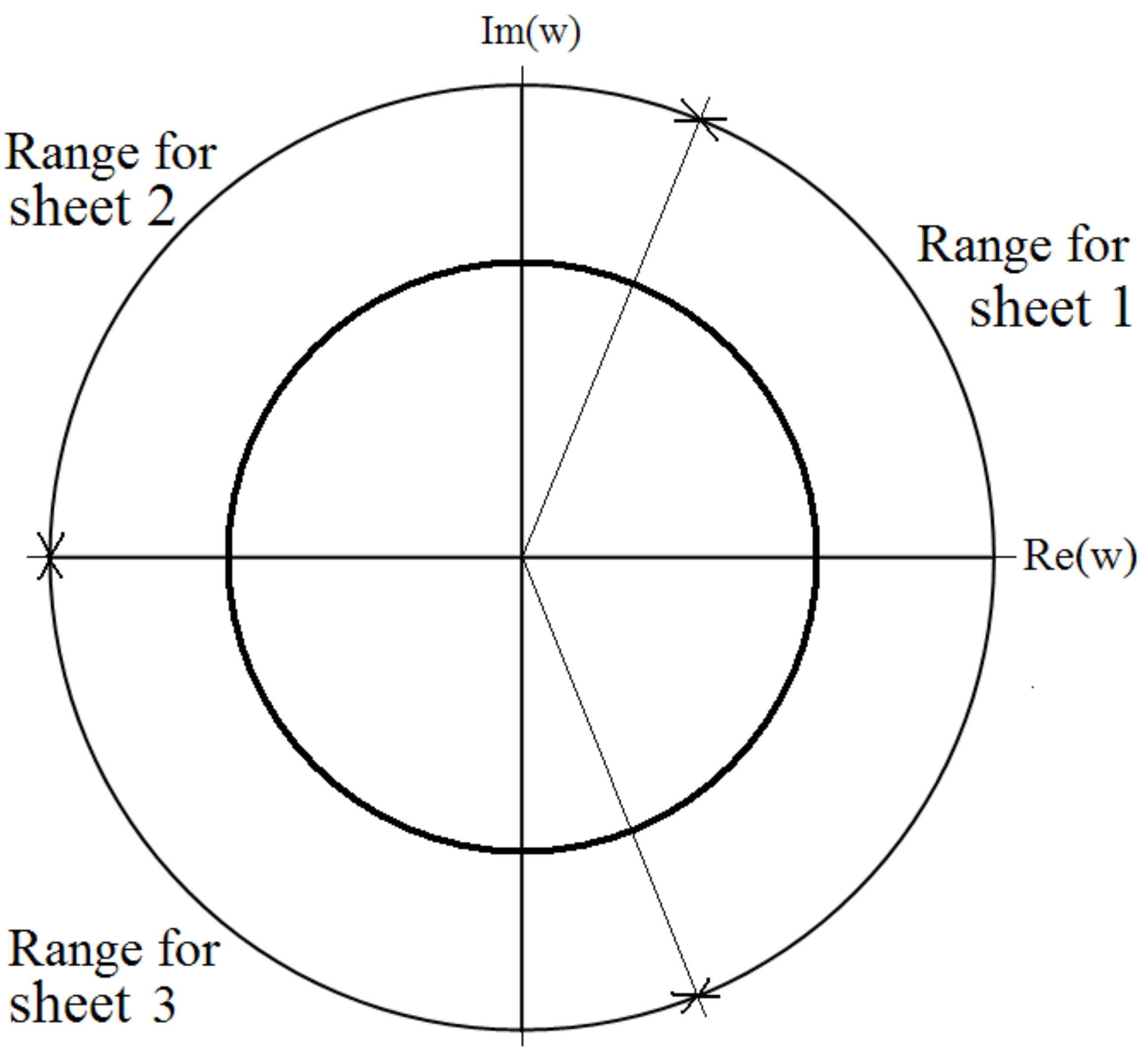}}
   \caption{Correspondence between the 3-sheets Riemann surface and $w$-plane for Eq.(\ref{w3}).}
   \label{Riemann3}
\end{figure}

\noindent
\textbf{Example 3}:
Analogous to previous examples we can also investigate function
\begin{equation}\label{w3}
  w=s^{\frac{1}{3}},
\end{equation}
where in this case the Riemann surface has three sheets and each maps onto one-third
of the $w$-plane (see Fig.~\ref{Riemann3}).

\begin{figure}[!h]
\centering
\noindent
  \subfigure[Riemann surface]{
  	\includegraphics[width=3.5in]{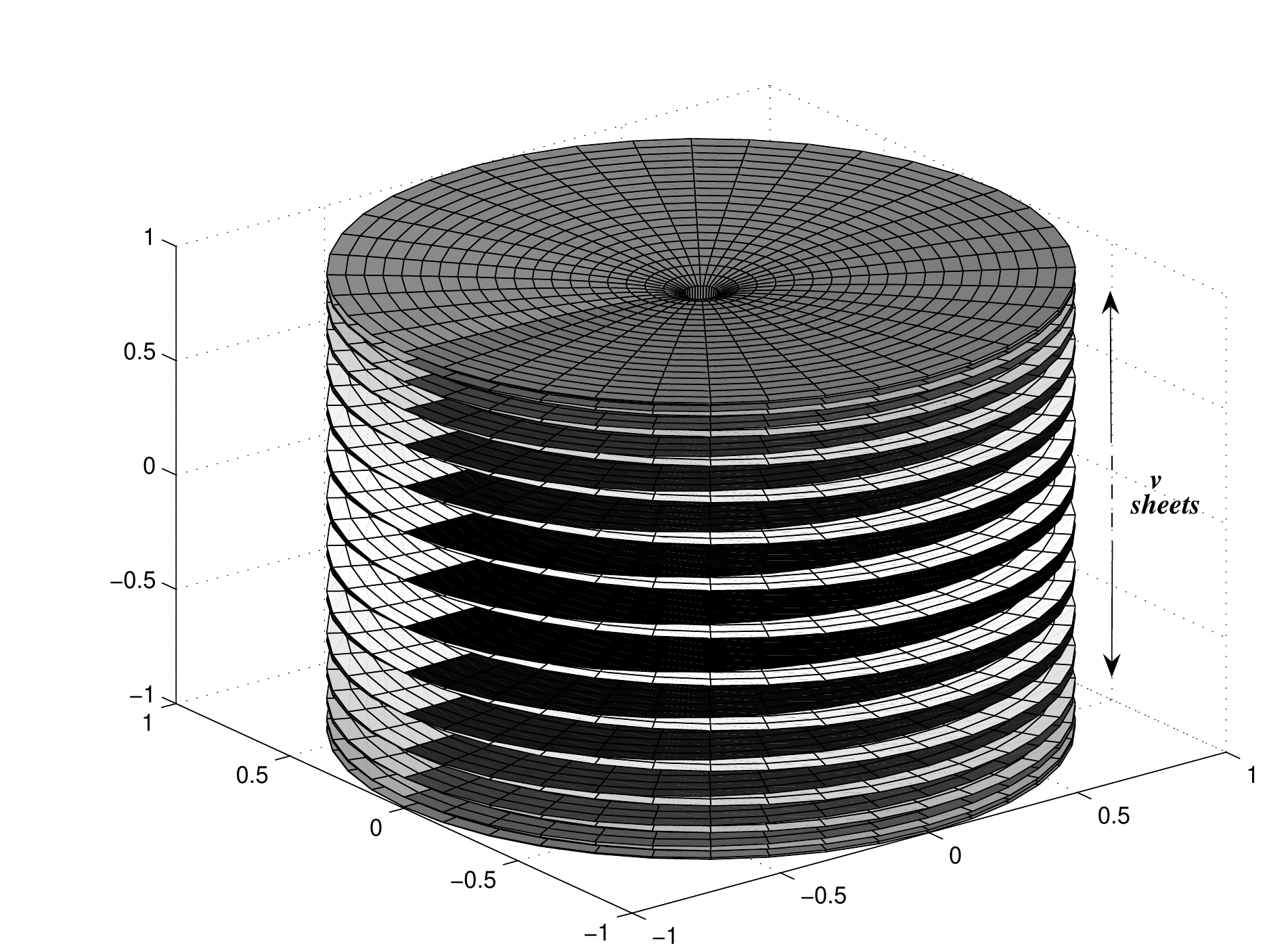}}
  \subfigure[Complex $w$-plane]{\label{RSn}
  	\includegraphics[width=2.9in]{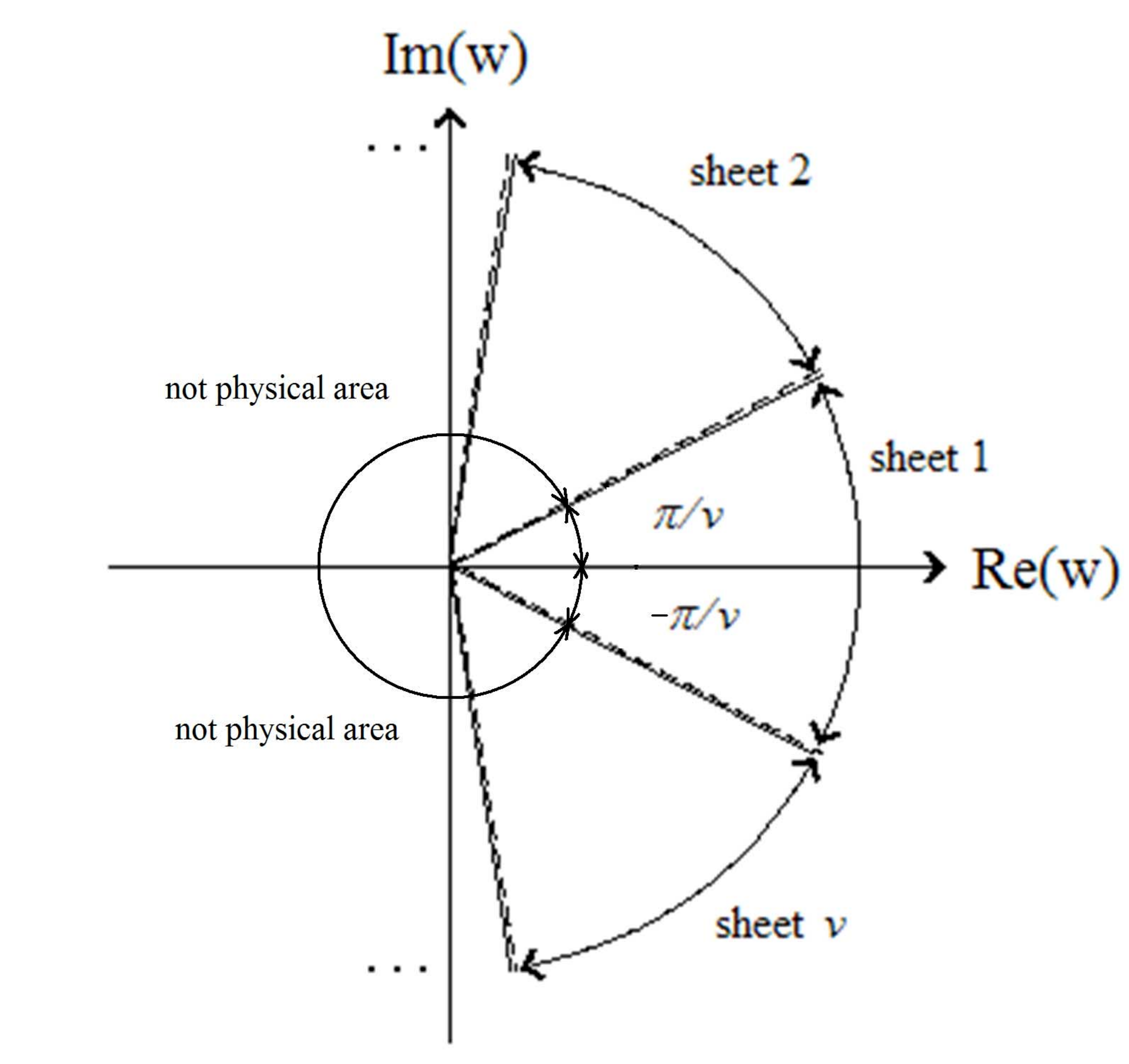}}
   \caption{Correspondence between the $w$-plane and the Riemann sheets for Eq.(\ref{v-sheets}).}
   \label{Riemann}
\end{figure}

\noindent
\textbf{Definition 1.} Generally, for the multivalued function defined as follow
\begin{equation}\label{v-sheets}
  w=s^{\frac{1}{v}},
\end{equation}
where $v \in N$ $(v=1,2,3,\dots)$ we get the $v$ sheets in the Riemann surface.
In Fig.\ref{Riemann} is shown the relationship between the $w$-plane and the $v$ sheets of the Riemann surface
where sector $-\pi/v < \mbox{arg} (w) \leq \pi/v$ corresponds to $\Omega$ (first Riemann sheet).

\begin{figure}[h]
\centering
\noindent
  \subfigure[$0 < q < 1$]{\label {s1}
		\includegraphics[width=2in]{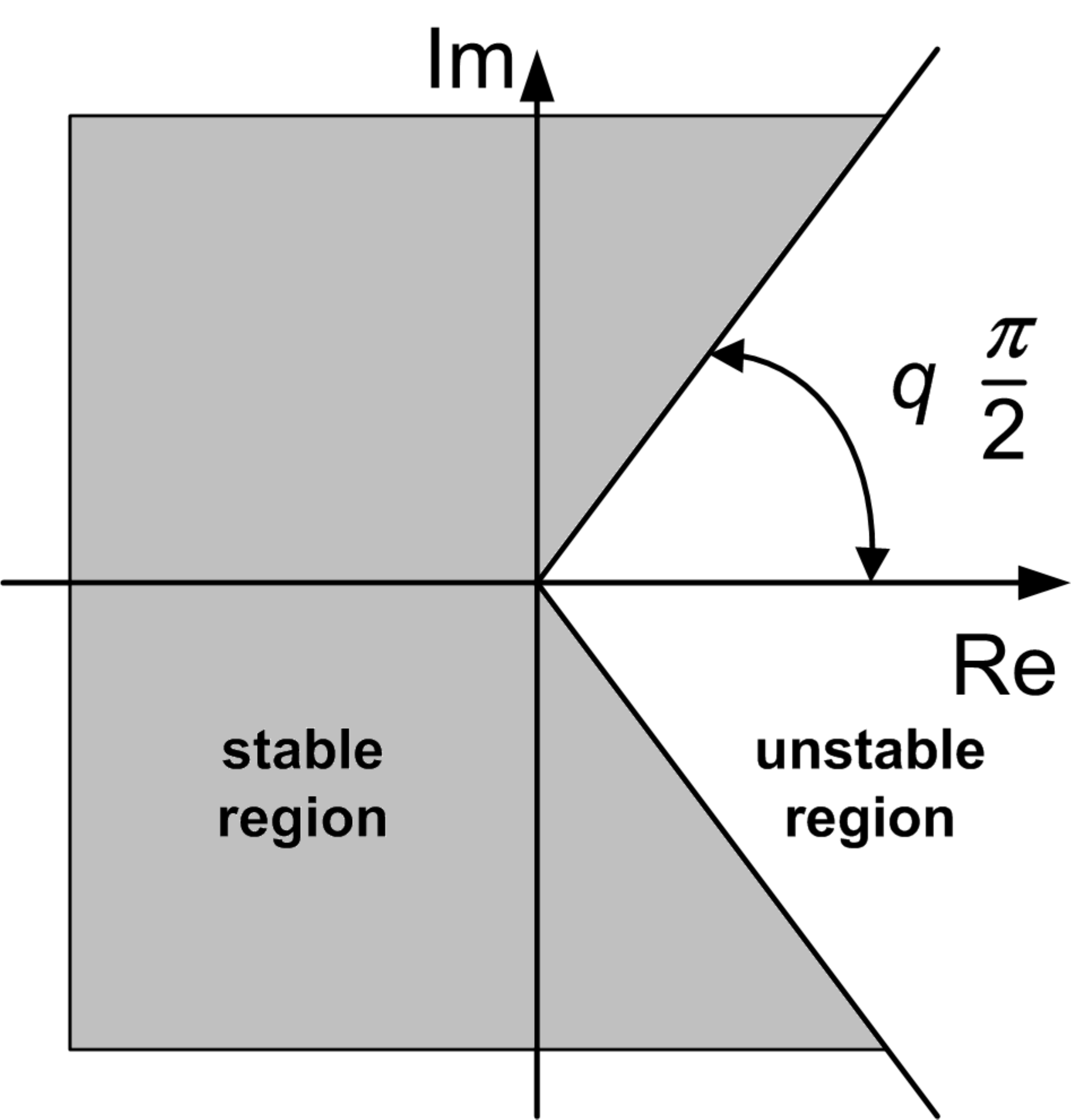}}
  \subfigure[$q = 1$]{\label {s2}
                \includegraphics[width=2in]{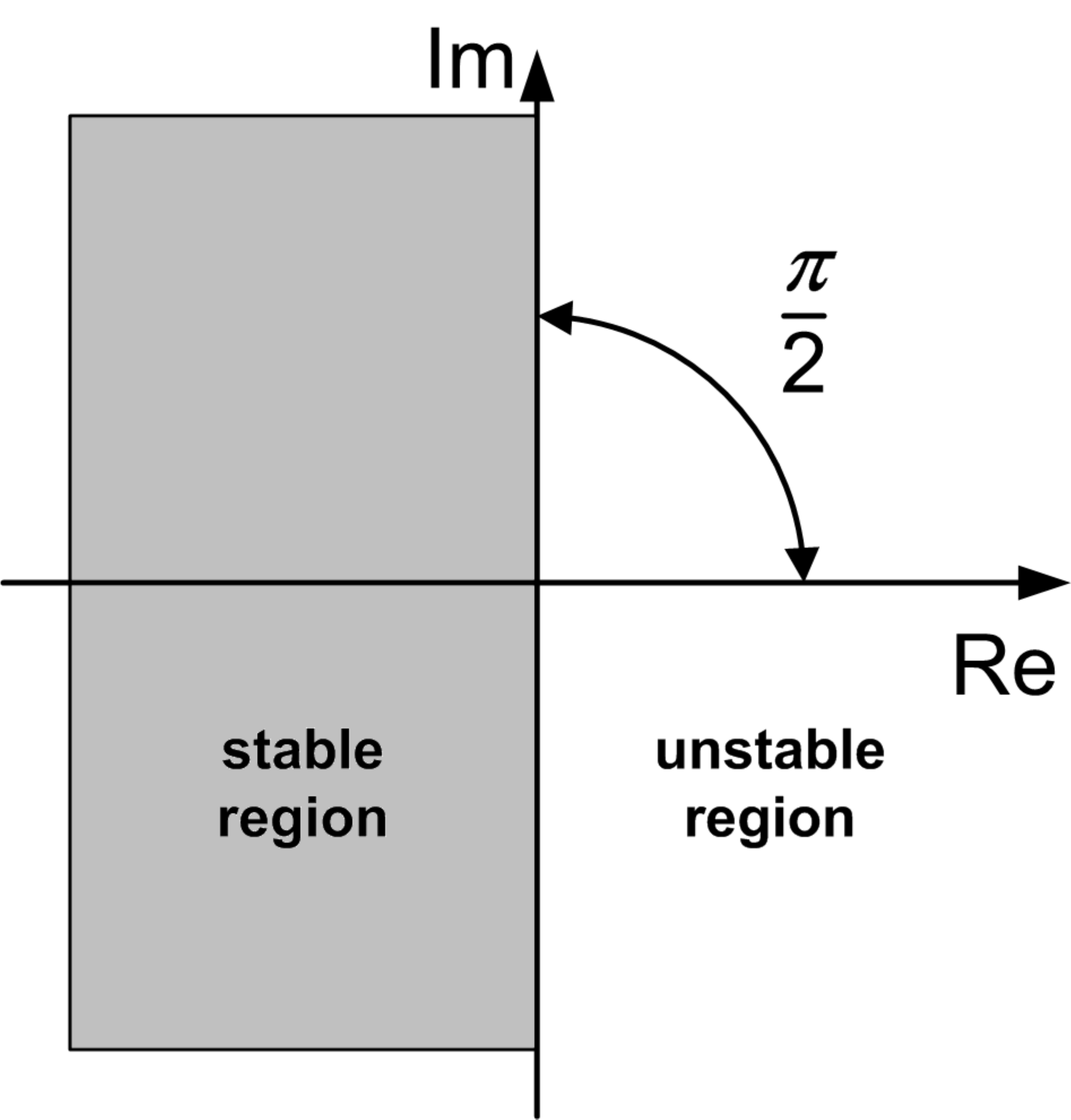}}
  \subfigure[$1 < q < 2$]{ \label {s3}
                \includegraphics[width=2in]{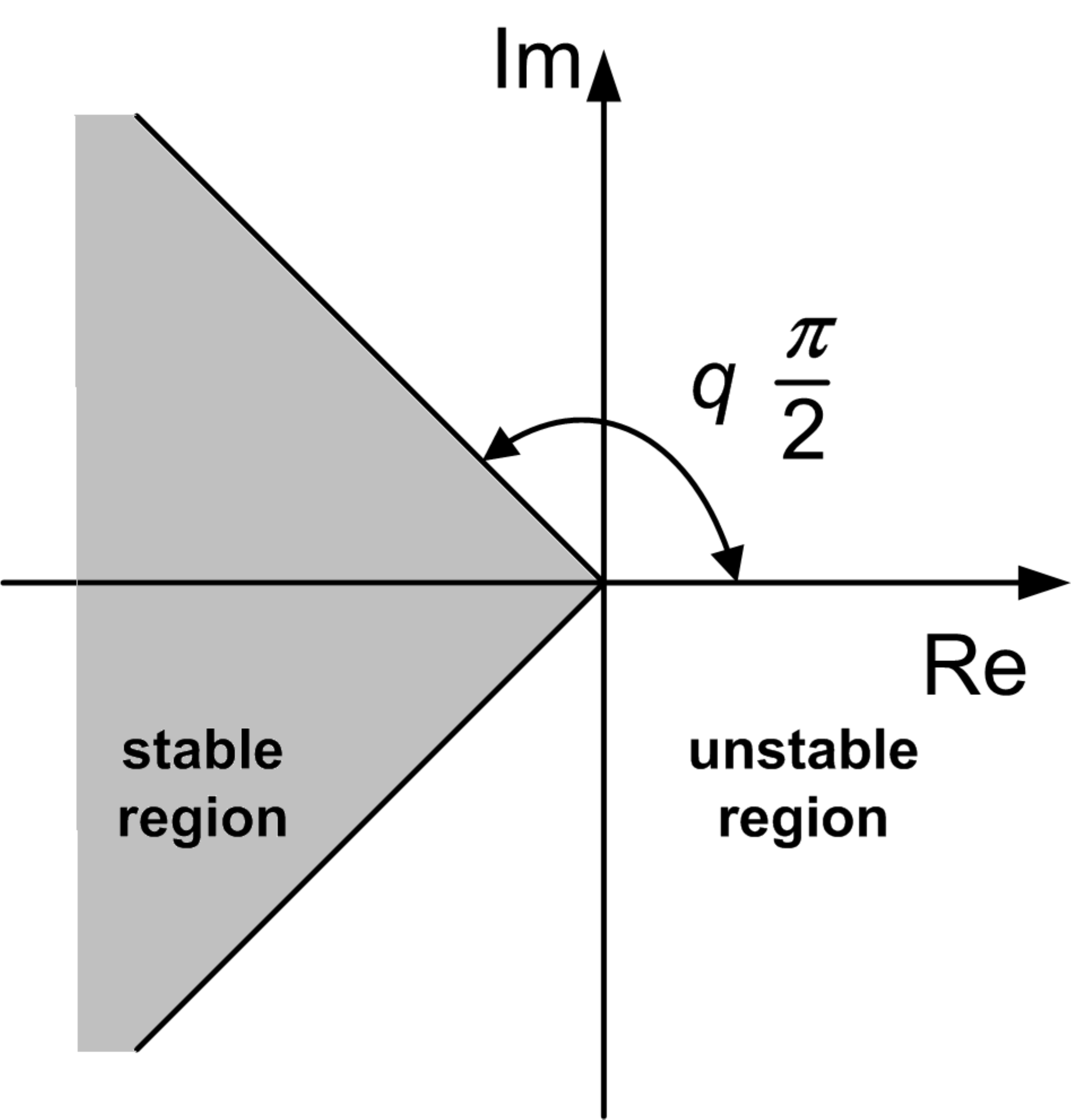}}
   \caption{Stability regions of the fractional order system.}
   \label{stability}
\end{figure}

\noindent
\textbf{Definition 2.} Mapping the poles from $s^{q}$-plane into the $w$-plane, where $q \in Q$ such as  $q = \frac{k}{m}$
for $k, m \in N$ and $|\mbox{arg}(w)|=|\phi|$, can be done by the following rule: If we assume $k=1$, then the mapping from $s$-plane
to $w$-plane is independent of $k$.
Unstable region from $s$-plane transforms to sector $|\phi|<\frac{\pi}{2m}$ and stable region transforms
to sector $\frac{\pi}{2m}<|\phi|<\frac{\pi}{m}$. The region where $|\phi|>\frac{\pi}{m}$ is not physical.
Therefore, the system will be stable if all roots in the $w$-plane lie in the region $|\phi|>\frac{\pi}{2m}$.
Stability regions depicted in Fig.~\ref{stability} correspond to the following propositions:

\begin{enumerate}
\item For $k<m$\, $(q<1)$ the stability region is depicted in Fig.~\ref{s1}.

\item For $k=m$\, $(q=1)$ the stability region corresponds to the $s$-plane (see Fig.~\ref{s2}).

\item For $k>m$\, $(q>1)$ the stability region is depicted in Fig.~\ref{s3}.
\end{enumerate}

\subsection{Stability of Fractional LTI Systems}

As we can see in previous subsection, in the fractional case, the stability is different from the integer one.
Interesting notion is that a stable fractional system may have roots in right half of complex $w$-plane (see Fig.~\ref{stability}).

Since the principal sheet of the Riemann surface is defined $-\pi < \mbox{arg}(s) < \pi$,
by using the mapping $w = s^q$, the corresponding $w$ domain is defined by $-q\pi < \mbox{arg}(w) < q\pi$,
and the $w$ plane region corresponding to the right half plane of this sheet is defined
by $-q\pi/2 < \mbox{arg}(w) < q\pi/2$.

Consider the fractional order pseudo-polynomial
$$
 Q(s) = a_1 s^{q_1} + a_2 s^{q_2} + \ldots + a_n s^{q_n}=
       a_1 s^{c_1/d_1} + a_2 s^{c_2/d_2} + \ldots + a_n s^{c_n/d_n},
$$
where $q_i$ are rational number expressed as $c_i/d_i$
and $a_i$ are the real numbers for $i=1,2,\dots, n$.
If for some $i$, $c_i=0$ then $d_i=1$. Let $v$ be the least common multiple (LCM) of
$d_1, d_2, \ldots d_n$ denote as $v=\mbox{LCM}\{d_1, d_2, \ldots d_n\}$, then \cite{Ghartemani}
\begin{equation} \label{eqLCM}
      Q(s) = a_1 s^{\frac{v_1}{v}} + a_2 s^{\frac{v2}{v}} + \ldots + a_n s^{\frac{v_n}{v}} =
             a_1 (s^{\frac{1}{v}})^{v_1} + a_2 (s^{\frac{1}{v}})^{v_2} + \ldots + a_n (s^{\frac{1}{v}})^{v_n}.
\end{equation}
The fractional degree (FDEG) of the polynomial $Q(s)$ is defined as \cite{Ghartemani}
$$
\mbox{FDEG}\{Q(s)\} = \mbox{max} \{ v_1, v_2, \dots, v_n\}.
$$
The domain of definition for (\ref{eqLCM}) is the Riemann surface with $v$ Riemann sheets where
origin is a branch point of order $v-1$ and the branch cut is assumed at $\mbox{R}^{-}$.
Number of roots for fractional algebraic equation (\ref{eqLCM}) is given by the following proposition \cite{Bayat1}: \\
\textbf{Proposition 1.} Let $Q(s)$ be a fractional order polynomial with $\mbox{FDEG}\{Q(s)\}=n$. Then the equation Q(s)=0 has
exactly $n$ roots on the Riemann surface \cite{Bayat1}. 

\noindent
\textbf{Definition 3.}  The fractional order polynomial
$$
  Q(s)=a_1 s^{\frac{n}{v}} + a_2 s^{\frac{n-1}{v}} + \ldots + a_n s^{\frac{1}{v}}  + a_{n+1}
$$
is \textit{minimal} if $\mbox{FDEG}\{Q(s)\}=n$. We will assume that all fractional order
polynomial are minimal.
This ensures that there is no redundancy
in the number of the Riemann sheets \cite{Ghartemani}.

On the other hand, it has been shown, by several authors and by using several methods,
that for the case of FOLTI system of commensurate order, a geometrical method of complex analysis based on
the argument principle of the roots of the characteristic equation (a polynomial in
this particular case) can be used for the stability check in the BIBO
sense (see e.g. \cite{Matignon, Petras2}).
The stability condition can then be stated as follows
\cite{Matignon96, Matignon, Vinagre5}:        \\
\textbf{Theorem 1.}
A commensurate order system described by a rational transfer
function (\ref{Eq2}) is stable if only if
$$\left | \arg \left( \lambda _{i}\right) \right | > \alpha \frac{\pi }{2},\,\,\mbox{for all}\,\, i$$
with $\lambda _{i}$ the $i$-th root of $P(s^{\alpha})$.

For the FOLTI system with commensurate order where the system poles are
in general complex conjugate, the stability condition can also be
expressed as follows \cite{Matignon96, Matignon}: \\
\textbf{Theorem 2.}
A commensurate order system described by a rational transfer
function
$$
G(w )=\frac{Q(w )}{P(w)},$$
where $w =s^{q},\, q \in \mbox{R}^{+}, (0<q<2)$, is stable if only if
$$\left| \arg \left( w _{i}\right) \right| >q \frac{\pi}{2},$$
with $\forall w_{i} \in \mbox{C}$ the $i$-th root of $P(w)=0$.

When $w=0$ is a single root (singularity at the origin) of $P$, the system cannot be stable. For $q=1$, this is the
classical theorem of pole location in the complex plane: have no pole in the closed right half plane of the first Riemann sheet.
The stability region suggested by this theorem tends to the whole $s$-plane when $q$ tends to 0, corresponds to the
Routh-Hurwitz stability when $q$ = 1, and tends to the negative real axis when $q$ tends to 2.

\noindent
\textbf{Theorem 3.}
It has been shown that commensurate system (\ref{LTI_SS}) is stable if the following condition is satisfied
(also if the triplet \textbf{A}, \textbf{B}, \textbf{C} is minimal) \cite{Aoun1, Matignon, Tavazoei, Haeri, Tavazoei2}:
\begin{equation}
|\mbox{arg(eig(\textbf{A}))}| > q \frac{\pi}{2},
\end{equation}
where $0<q<2$ and eig(\textbf{A}) represents the eigenvalues of matrix \textbf{A}. \\

\noindent
\textbf{Proposition 2.}
We can assume, that some incommensurate order systems described by the FODE (\ref{n-term-equation_u})
or (\ref{LTI_SS}), %
can be decomposed to the following modal form of the fractional transfer function
(so called Laguerre functions \cite{Aoun2}):
\begin{equation}\label{ds1}
F(s) = \sum_{i=1}^{N} \sum_{k=1}^{n_k} \frac{A_{i,k}}{(s^{q_i}+\lambda_i)^k}
\end{equation}
for some complex numbers $A_{i,k}$, $\lambda_i$, and positive integer $n_k$.

A system (\ref{ds1}) is BIBO stable if and only if $q_i$
and  the argument of $\lambda_i$ denoted by arg($\lambda_i$) in (\ref{ds1})
satisfy the inequalities
\begin{equation}\label{genstability}
0<q_i<2\,\,\,\,\,\mbox{and}\,\,\,\, \left | \arg \left( \lambda_i \right) \right | <
 \pi \left( 1- \frac{q_i}{2}\right)\,\,\,\,\ \mbox{for all}\,\, i.
\end{equation}
Henceforth, we will restrict the parameters $q_i$ to the interval $q_i \in (0, 2)$.  For the case $q_i=1$ for all $i$ we obtain
a classical stability condition for integer order system (no pole is in right half plane).
The inequalities (\ref{genstability}) were obtained by applying the stability
results given in \cite{Akcay, Matignon}.
\noindent
\textbf{Theorem 4.}
Consider the following autonomous system for internal stability definition \cite{Deng}:
\begin{equation}\label{LTI_SS_2}
_0D^{\textbf{q}}_{t} x(t) = \textbf{A} x(t), \,\,\,\, x(0)=x_0,
\end{equation}
with $\textbf{q}$ = $[q_1, q_2, \dots, q_n]^T$ and its $n$-dimensional representation:
\begin{eqnarray} \label{LTI_SS_3}
_0D^{q_1}_{t} x_1(t) &=& a_{11}x_1(t) + a_{12}x_2(t) + \dots + a_{1n}x_n(t) \nonumber \\
_0D^{q_2}_{t} x_2(t) &=& a_{21}x_1(t) + a_{22}x_2(t) + \dots + a_{2n}x_n(t) \nonumber \\
\dots \nonumber \\
_0D^{q_n}_{t} x_n(t) &=& a_{n1}x_1(t) + a_{n2}x_2(t) + \dots + a_{nn}x_n(t)
\end{eqnarray}
where all $q_i$'s are rational numbers between $0$ and $2$. Assume $m$ be the LCM
of the denominators $u_i$'s of $q_i$'s , where $q_i=v_i/u_i$, $v_i, u_i \in Z^+$
for $i=1,2, \dots, n$ and we set $\gamma = 1/m$. Define:
\begin{equation}\label{chaeq}
\mbox{det}
\left(
\begin{array}{cccc}
\lambda^{mq_1} - a_{11} & -a_{12} & \dots & -a_{1a}  \\
-a_{21} & \lambda^{mq_2} - a_{22} & \dots &  -a_{2n} \\
\dots       \\
-a_{n1} & -a_{n2} & \dots & \lambda^{mq_n} - a_{nn}
\end{array}
\right)  =0.
\end{equation}
The characteristic equation (\ref{chaeq}) can be transformed to integer order polynomial equation
if all $q_i$'s are rational number. Then the zero solution of system  (\ref{LTI_SS_3}) is globally asymptotically
stable if all roots $\lambda_i$'s of the characteristic (polynomial) equation (\ref{chaeq}) satisfy
$$
|\mbox{arg}(\lambda_i)| > \gamma \frac{\pi}{2}\,\,\mbox{for all}\,\, i.
$$
Denote $\lambda$ by $s^{\gamma}$ in equation (\ref{chaeq}), we get
the characteristic equation in the form $\mbox{det}(s^{\gamma}I-A)=0$ and this assumption
was proved in paper \cite{Deng}.
\\
\textbf{Corollary 1.} Suppose $q_1=q_2=\dots, q_n \equiv q$, $q\in(0,2)$,
all eigenvalues $\lambda$ of matrix $A$ in (\ref{LTI_SS_M}) satisfy $|\mbox{arg}(\lambda)| >q \pi /2$,
the characteristic equation becomes $\mbox{det}(s^qI-A)=0$ and all
characteristic roots of the system (\ref{LTI_SS}) have negative real parts \cite{Deng}.
This result is Theorem 1 of paper \cite{Matignon96}.
\\
\noindent
\textbf{Remark 1.} Generally, when we assume $s=|r|e^{i \phi}$, where $|r|$ is modulus and $\phi$
is argument of complex number in $s$-plane, respectively, transformation $w=s^{\frac{1}{m}}$
to complex $w$-plane can be viewed as  $s=|r|^{\frac{1}{m}}e^{\frac{i \phi}{m}}$ and thus
$|\mbox{arg}(s)|$ = $m. |\mbox{arg}(w)|$ and $|s|=|w|^m$. Proof of this statement is obvious.
\\
\noindent
\textbf{Stability analysis criteria for a general FOLTI system can be summarized as follow:}
\\
The characteristic equation of a general LTI fractional order system of the form:
\begin{equation}
a_{n}s^{\alpha_{n}}+\ldots +a_{1}s^{\alpha_{1}}+a_{0}s^{\alpha_{0}} \equiv \sum_{i=0}^{n} a_i s^{\alpha_i} = 0
\end{equation}
may be rewritten as
$$
\sum_{i=0}^{n} a_i s^{\frac{u_i}{v_i}} = 0
$$
and transformed into $w$-plane
\begin{equation}  \label{cheqw}
\sum_{i=0}^{n} a_i w^i = 0,
\end{equation}
with $w=s^{\frac{k}{m}}$, where $m$ is the LCM of $v_i$.
The procedure of stability analysis is (see e.g. \cite{Radwan}):
\begin{enumerate}
\item For given $a_i$ calculate the roots of Eq.(\ref{cheqw}) and find the absolute phase of all roots $|\phi_w|$.
\item Roots in the primary sheet of the $w$-plane which have corresponding roots in the $s$-plane can
be obtained by finding all roots which lie in the region $|\phi_w|<\frac{\pi}{m}$ then applying
the inverse transformation $s=w^m$ (see Remark 1.). The region where $|\phi_w|>\frac{\pi}{m}$ is not physical.
For testing the roots in desired region the matrix approach can be used (\ref{matrixstability}).
\item The condition for stability is $\frac{\pi}{2 m}<|\phi_w|<\frac{\pi}{m}$. Condition for
oscillation is $|\phi_w|=\frac{\pi}{2 m}$ otherwise the system is unstable (see Fig.~\ref{RSn}).
If there is not root in the physical $s$-plane, the system will always be stable \cite{Radwan}.
\end{enumerate}

\noindent
\textbf{Example 4.}
Let us consider the linear fractional order LTI system described by the transfer function
\cite{Dorcak, Podlubny}:
\begin{equation}
G(s)=\frac{Y(s)}{U(s)}=\frac{1}{0.8 s^{2.2} + 0.5 s^{0.9} +1},
\end{equation}
and corresponding FODE has the following form:
\begin{equation} \label{Ex1}
0.8\, _0D^{2.2}_t y(t) + 0.5\,_0D^{0.9}_t y(t) + y(t) = u(t)
\end{equation}
with zero initial conditions.

The system (\ref{Ex1}) can be rewritten to its state space representation $(x_1(t) \equiv y(t))$:
\begin{eqnarray} \label{LTI_SS_Ex2}
\left[
\begin{array}{c}
_0D^{\frac{9}{10}}x_{1}(t) \\
_0D^{\frac{13}{10}}x_{2}(t) \\
\end{array}
\right] &=&\left[
\begin{array}{cc}
0 & 1  \\
-1/0.8 & -0.5/0.8
\end{array}
\right] \left[
\begin{array}{c}
x_{1}(t) \\
x_{2}(t) \\
\end{array}
\right] +\left[
\begin{array}{c}
0 \\
1/0.8
\end{array}
\right] u(t) \nonumber \\
y(t) &=&\left[
\begin{array}{cc}
1 & 0
\end{array}
\right] \left[
\begin{array}{c}
x_{1}(t) \\
x_{2}(t) \\
\end{array}
\right]
\end{eqnarray}
The eigenvalues of the matrix \textbf{A} are $\lambda_{1,2} = -0.3125 \pm 1.0735j$
and then $|\mbox{arg}(\lambda_{1,2})| = 1.8541$. Because of various derivative
orders in (\ref{LTI_SS_Ex2}), the Theorem 3 cannot be used directly.

\begin{figure}[ht]
\centering
\noindent
  \includegraphics[width=4.4in]{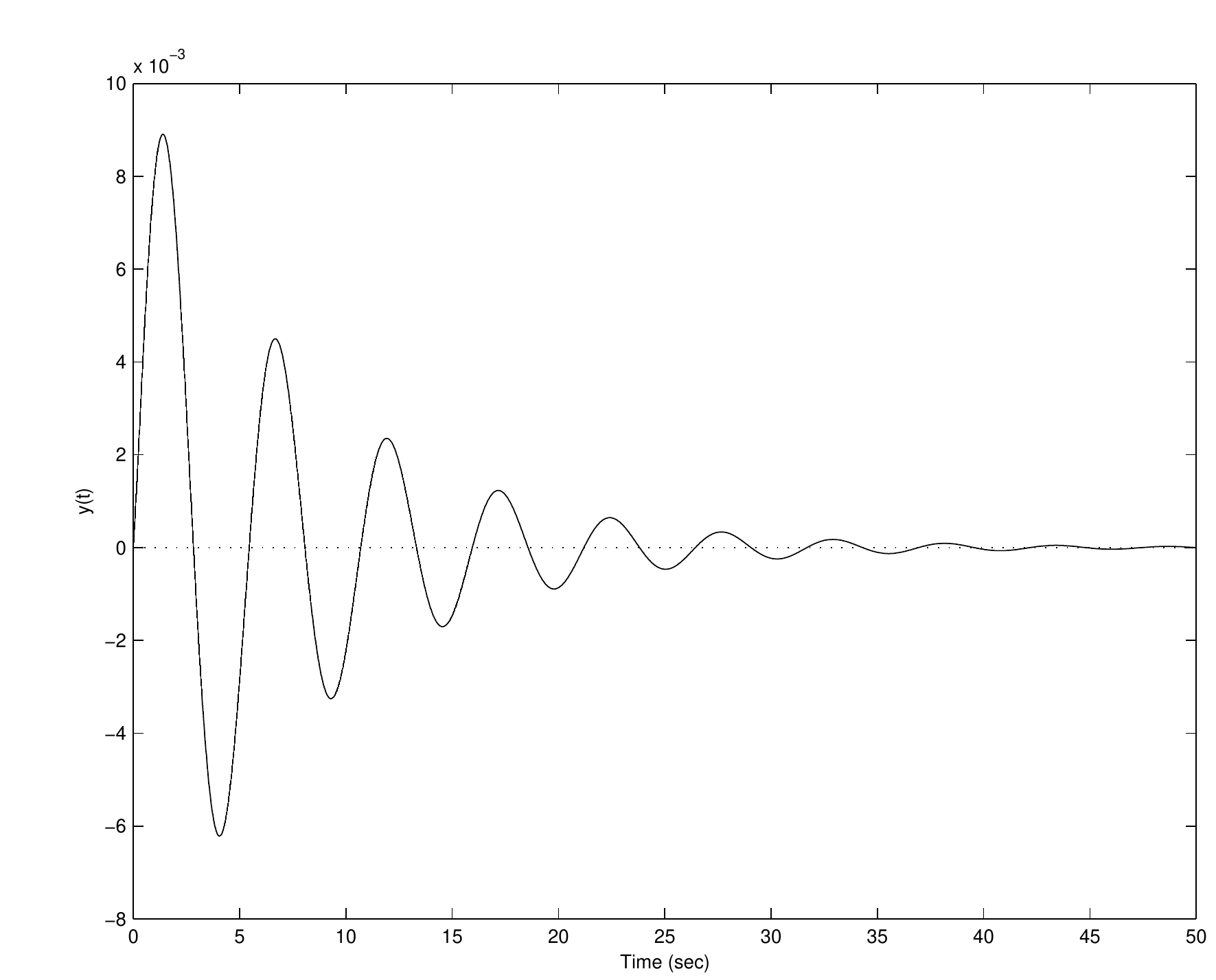}
   \caption{Analytical solution of the FODE (\ref{Ex1}) where $u(t)=0$ for $50\,sec$.}
   \label{impulse1}
\end{figure}

The analytical solution of the FODE (\ref{Ex1}) for $u(t)=0$ obtained from general solution (\ref{Step-response}) has form:
\begin{equation}
y(t) = \frac{1}{0.8} \sum_{k=0}^{\infty} \frac{(-1)^k}{k !} \left ( \frac{1}{0.8} \right )^k {\cal E}_{k}(t, -\frac{0.5}{0.8}; 2.2-0.9, 2.2 + 0.9k).
\end{equation}

In Fig.~\ref{impulse1} is depicted the analytical solution of the FODE (\ref{Ex1}) where $u(t)=0$.
As we can see in the figure, solution is stable because $\mbox{lim}_{t \to \infty} y(t) =0$.
Let us investigate stability according to the previously described method.
The corresponding characteristic equation of system is:
\begin{equation}\label{ChEq1}
P(s): 0.8 s^{2.2} + 0.5 s^{0.9} +1 = 0 \,\,\,\Rightarrow \,\,\, 0.8 s^{\frac{22}{10}} + 0.5 s^{\frac{9}{10}} +1 = 0,
\end{equation}
when $m=10$, $w=s^{\frac{1}{10}}$ then the roots $w_i$'s and their appropriate arguments of polynomial
\begin{equation}\label{roots10}
P(w): 0.8 w^{22} + 0.5 w^{9} + 1 = 0
\end{equation}
are:    \\ %
$
w_{1,2}=-0.9970 \pm 0.1182j, |\mbox{arg}(w_{1,2})|= 3.023; w_{3,4}=-0.9297 \pm 0.4414j, |\mbox{arg}(w_{3,4})|=2.698; \\
w_{5,6}=-0.7465 \pm 0.6420j, |\mbox{arg}(w_{5,6})|= 2.431; w_{7,8}=-0.5661 \pm 0.8633j, |\mbox{arg}(w_{7,8})|=2.151; \\
w_{9,10}=-0.259 \pm 0.9625j, |\mbox{arg}(w_{9,10})|=1.834; w_{11,12}= -0.0254 \pm 1.0111j, |\mbox{arg}(w_{11,12})|=1.595; \\
w_{13,14}=0.3080 \pm 0.9772j, |\mbox{arg}(w_{11,12})|=1.265; w_{15,16}=  0.5243 \pm 0.8359j, |\mbox{arg}(w_{15,16})|=1.010; \\
w_{17,18}=0.7793 \pm 0.6795j, |\mbox{arg}(w_{17,18})|= 0.717;  w_{19,20}=0.9084 \pm 0.3960j, |\mbox{arg}(w_{19,20})|=0.411;    \\
w_{21,22}=1.0045 \pm 0.1684j,  |\mbox{arg}(w_{21,22})| = 0.1661;  \\
$
\\
Physical significance roots are in the first Riemann sheet, which is expressed by
relation $-\pi/m < \phi <\pi/m$, where $\phi = \mbox{arg}(w)$. In this case
they are complex conjugate roots $w_{21,22}=1.0045 \pm 0.1684j$\, ($|\mbox{arg}(w_{21,22})|=0.1661$), which satisfy
conditions $|\mbox{arg}(w_{21,22})|>\pi/2m=\pi/20$. It means that system (\ref{Ex1}) is stable (see Fig.~\ref{Riemann10}).
Other roots of the polynomial equation (\ref{roots10}) lie in region $|\phi|>\frac{\pi}{m}$ which is not physical
(outside of closed angular sector limited by thick line in Fig.~{\ref{Riemann10b}).
\vspace{-3mm}
\begin{figure}[ht]
\centering
\noindent
  \subfigure[10-sheets Riemann surface]{ \label{Riemann10a}
  \includegraphics[width=3.3in]{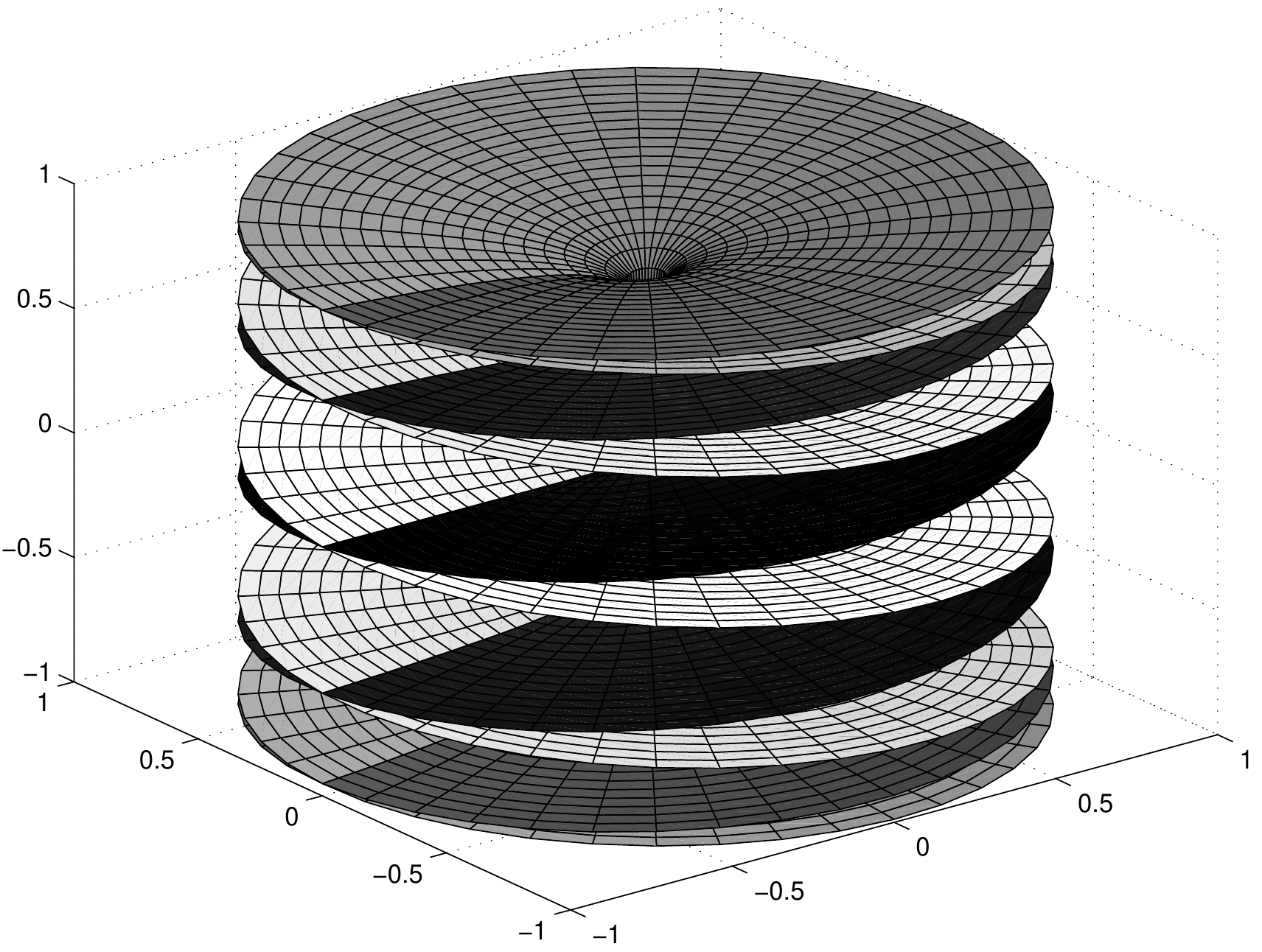}}
  \subfigure[Poles in complex $w$-plane]{ \label{Riemann10b}
  \includegraphics[width=2.9in]{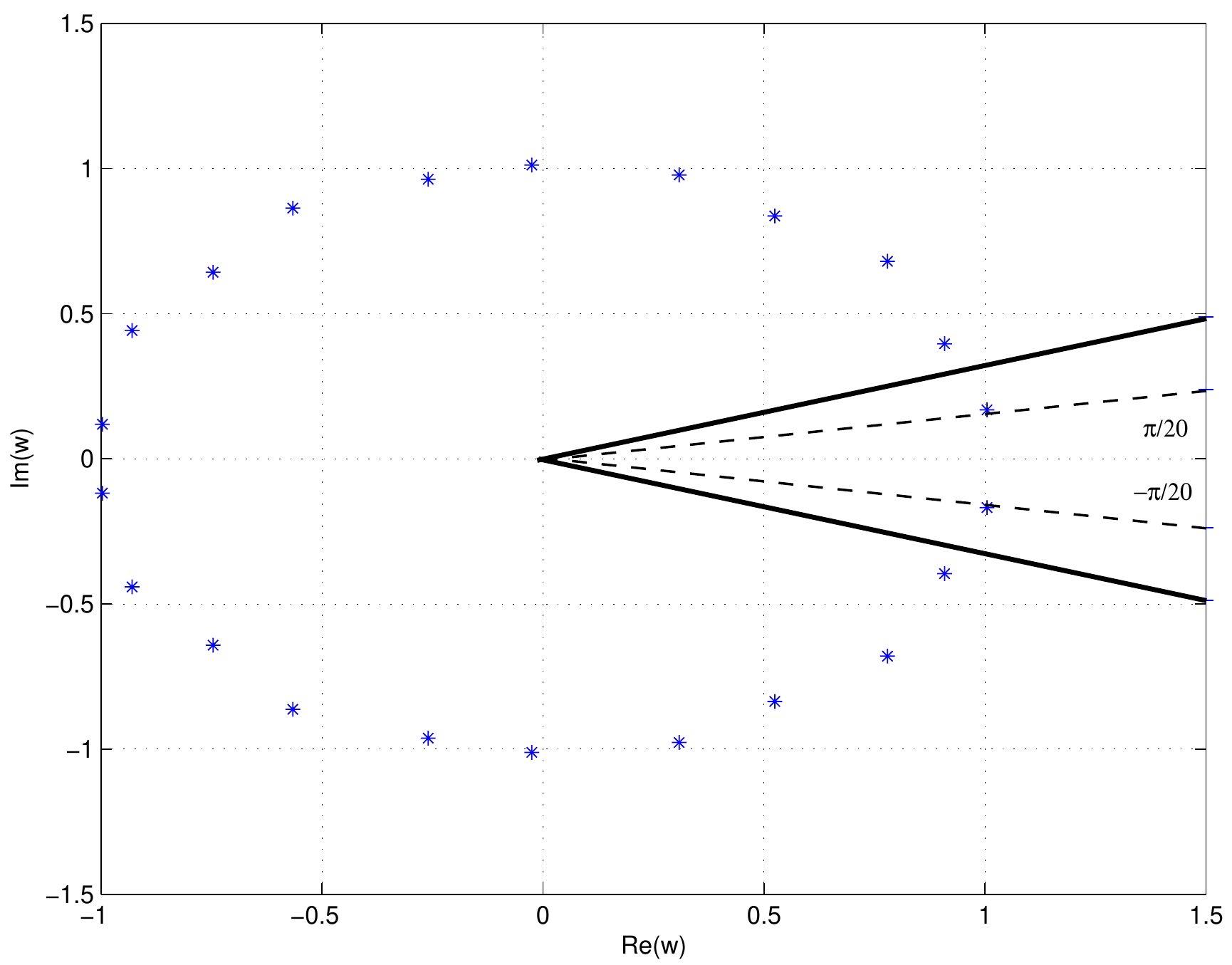}}
  \caption{Riemann surface of function $w=s^{\frac{1}{10}}$ and roots of Eq.(\ref{roots10}) in complex $w$-plane.}
   \label{Riemann10}
\end{figure}

In Fig.~\ref{Riemann10a} is depicted the Riemann surface of the function $w=s^{\frac{1}{10}}$ with the 10-Riemann sheets
and in Fig.~\ref{Riemann10b} are depicted the roots in complex $w$-plane with angular sector corresponds to stability
region (dashed line) and the first Riemann sheet (thick line).

The interesting notion of Remark 1 should be mentioned here. The characteristic equation (\ref{ChEq1}) has the following poles:
$$
  s_{1,2} = -0.10841 \pm 1.19699j,
$$
in the first Riemann sheet in $s$-plane, which can be obtained e.g. via the Matlab routine as for instance:
\begin{verbatim}
>>s=solve('0.8*s^2.2+0.5*s^0.9+1=0','s')
\end{verbatim}
When we compare $|\mbox{arg}(w_{21,22})|=0.1661$ and $|\mbox{arg}(s_{1,2})|=1.661$, we can see
that $|\mbox{arg}(s_{1,2})| = m |\mbox{arg}(w_{21,22})|$, where $m=10$ in transformation
$w=s^{\frac{1}{m}}$. The first Riemann sheet is transformed from $s$-plane to $w$-plane
as follow: $-\pi/10 < \mbox{arg}(w) < \pi/10$ and in order to $-\pi < 10 . \mbox{arg}(w) < \pi$.
Therefore from this consideration we then obtain $|\mbox{arg}(s)|$ = $10. |\mbox{arg}(w)|$.

\noindent
\textbf{Example 5.}
Let us examine an interesting example of application, so called Bessel function of the
first kind, which transfer function is \cite{Matignon}:
\begin{equation}\label{bassel}
      H(s) = \frac{1}{\sqrt{s^2 + 1}}\,\,\,\,\, \forall s, \,\, \Re(s) >0.
\end{equation}
We have two branch points $s_1=i$, and $s_2=-i$ and two cuts. One along the half line
$(-\infty +i, i)$ and another one along the half line $(-\infty - i, -i)$. In this
doubly cut complex plane, we have the identity $\sqrt{s^2 +1} = \sqrt{s-i}\sqrt{s+i}$.
The well known asymptotic expansion of Eq.(\ref{bassel}) is:
$$
  h(t) \approx \sqrt{\frac{2}{\pi t}} \mbox{cos}(t - \frac{\pi}{4}) = \sqrt{\frac{2}{\pi}}\, t^{-\frac{1}{2}}\, E_{2, 1} \left(-(t-\frac{\pi}{4})^2\right).
$$
According to the branch points and above asymptotic expansion we can
state, that system described by the Bessel function (\ref{bassel})
is on boundary of stability and has oscillation behaviour.

\noindent
\textbf{Example 6.}
Consider the closed loop system with controlled system (electrical heater)
\begin{equation}
G(s) = \frac{1}{39.96s^{1.25} + 0.598}
\end{equation}
and fractional order controller
\begin{equation}
C(s) = 64.47 + 12.46 s
\end{equation}
The resulting closed loop transfer function $G_c(s)$ becomes \cite{Petras5}:
\begin{equation}\label{G_s_u_c}
 G_c(s) = \frac{Y(s)}{W(s)}= \frac{12.46 s + 64.47}{39.69 s^{1.25}+12.46 s + 65.068}
\end{equation}
The analytical solution (impulse response) of the fractional order control system (\ref{G_s_u_c}) is:
\begin{eqnarray}
  y(t) & = & \frac{12.46}{39.69}\sum_{k=0}^{\infty} \frac{(-1)^k}{k!}\,
           \left( \frac{12.46}{39.69} \right)^k \,
\times 	{\cal E}_k(t, -\frac{65.068}{39.69}; 1.25, 0.25 - k)
		\nonumber \\
       & + & \frac{64.47}{39.69}\sum_{k=0}^{\infty} \frac{(-1)^k}{k!}\,
           \left ( \frac{65.068}{39.69} \right)^k \,
\times {\cal E}_k(t, -\frac{12.46}{39.69}; 1.25 - 1, 1.25 + k)
\end{eqnarray}
with zero initial conditions.

The characteristic equation of this system is
\begin{equation}
39.69 s^{1.25}+12.46 s + 65.068 = 0\,\,\, \Rightarrow\,\,\, 39.69 s^{\frac{5}{4}}+12.46 s^{\frac{4}{4}} + 65.068 = 0
\end{equation}
Using the notation $w=s^{\frac{1}{m}}$, where LCM is $m=4$, we obtain a polynomial of complex variable $w$ in form
\begin{equation}\label{poex6}
 39.69 w^5 + 12.46 w^4 + 65.068 = 0.
\end{equation}
Solving the polynomial (\ref{poex6}) we get the following roots and their arguments:
$$
w_1=-1.17474,  |\mbox{arg}(w_{1})|= \pi
$$
$$
w_{2,3} =  -0.40540 \pm 1.0426j,   |\mbox{arg}(w_{2,3})|= 1.9416
$$
$$
w_{4,5} = 0.83580 \pm 0.64536j,     |\mbox{arg}(w_{4,5})|=  0.6575
$$
This first Riemann sheet is defined as a sector in $w$-plane within interval
$
-\pi/4 < \mbox{arg}(w) < \pi/4.
$
Complex conjugate roots $w_{4,5}$ lie in this interval and satisfies the
stability condition given as $|\mbox{arg}(w)|>\frac{\pi}{8}$, therefore system is stable.
The region where $|\mbox{arg}(w)|>\frac{\pi}{4}$ is not physical.

\subsection{Stability of Fractional Nonlinear Systems}

As it was mentioned in \cite{Matignon96}, exponential stability cannot be used to characterize
asymptotic stability of fractional order systems. A new definition was introduced \cite{Oustaloup2}.

\noindent
\textbf{Definition 3.} Trajectory $x(t) = 0$ of the system  (\ref{FONS}) is $t^{-q}$ asymptotically
stable if there is  a~positive real $q$ such that:
$$
\forall ||x(t)||\,\,\, \mbox{with}\,\,\, t \leq t_0,\, \exists \, N(x(t)),\,\,\, \mbox{such that}\,\,\,\, \forall t \geq t_0, ||x(t)|| \leq N t^{-q}.
$$
The fact that the components of $x(t)$ slowly decay towards $0$ following $t^{-q}$ leads to fractional
systems sometimes being called long memory systems. Power law stability $t^{-q}$ is a special case
of the Mittag-Leffler stability \cite{Li}.

According to stability theorem defined in \cite{Tavazoei6},
the equilibrium points are asymptotically stable for $q_1 = q_2 = \dots = q_n \equiv q$ if all the eigenvalues $\lambda_i,\,(i=1,2,\dots,n)$
of the Jacobian matrix
$\textbf{J}=\partial \textbf{f} / \partial \textbf{x}$, where $\textbf{f}=[f_1,\, f_2,\ \dots,\ f_n]^T$, evaluated at the equilibrium,
satisfy the condition \cite{Haeri, Tavazoei}:
\begin{equation}\label{gen_stability}
|\mbox{arg(eig(\textbf{J}))}| = |\mbox{arg}(\lambda_i)| > q \frac{\pi}{2}, \,\,\,i=1,2,\dots,n.
\end{equation}
Fig.~\ref{stability} shows stable and unstable regions of the complex plane for such case.

Now, consider the incommensurate fractional order system $q_1 \neq q_2 \neq \dots \neq q_n$ and suppose
that $m$ is the LCM of the denominators $u_i$'s of $q_i$'s, where $q_i=v_i/u_i$, $v_i, u_i \in Z^+$
for $i=1,2, \dots, n$ and we set $\gamma = 1/m$. System (\ref{general}) is asymptotically stable if:
$$
|\mbox{arg}(\lambda)| > \gamma \frac{\pi}{2}
$$
for all roots $\lambda$ of the following equation
\begin{equation} \label{cheNS}
          \mbox{det}(\mbox{diag}([\lambda^{m q_1}\, \lambda^{m q_2}\,\ldots\,\lambda^{mq_n}])-\textbf{J} )=0.
\end{equation}

A necessary stability condition for fractional order systems (\ref{general}) to remain chaotic
is keeping at least one eigenvalue $\lambda$ in the unstable region \cite{Haeri}. The number of
saddle points and eigenvalues for one-scroll, double-scroll and multi-scroll attractors was
exactly described in work \cite{Tavazoei2}. Assume that 3D chaotic system has only three equilibria.
Therefore, if system has double-scroll attractor, it has two saddle points surrounded by scrolls
and one additional saddle point. Suppose that the unstable eigenvalues of scroll saddle points are:
$\lambda_{1,2} = \alpha_{1,2} \pm j \beta_{1,2}$.
The necessary condition to exhibit double-scroll attractor of system (\ref{general}) is the eigenvalues
$\lambda_{1,2}$ remaining in the unstable region \cite{Tavazoei2}. The condition for commensurate
derivatives order is
\begin{equation} \label{CO}
            q> \frac{2}{\pi} \mbox{atan} \left ( \frac{|\beta_i|}{\alpha_i} \right ), \,\,\,i=1,2.
\end{equation}
This condition can be used to determine the minimum order for which a nonlinear system can
generate chaos \cite{Haeri}.
\\
\textbf{Example 7.}
Let us investigate the Chen system with a double scroll attractor. The fractional order
form of such system can be described as \cite{Tavazoei6}
\begin{eqnarray}\label{ChenSystem}
   D^{0.8}_t x_1(t) &=& 35 [x_2(t) - x_1(t)] \nonumber \\
   D^{1.0}_t x_2(t) &=& - 7 x_1(t) -x_1(t) x_3 (t) + 28 x_2(t) \nonumber \\
   D^{0.9}_t x_3(t) &=& x_1(t) x_2(t) - 3 x_3(t)
\end{eqnarray}
The system has three equilibrium at $(0, 0, 0)$, $(7.94, 7.94, 21)$, and $(-7.94, -7.94, 21)$.
The Jacobian matrix of the system evaluated at $(x^{*}_1, x^{*}_2, x^{*}_3)$ is:
\begin{eqnarray}
\mbox{\textbf{J}} = \left[
\begin{array}{ccc}
-35 & 35 & 0 \\
-7-x^{*}_3 & 28 & -x^{*}_1 \\
x^{*}_2 & x^{*}_1 & -3
\end{array}
\right].
\end{eqnarray}
The two last equilibrium points are saddle points and surrounded by a chaotic double scroll attractor.
For these two points, equation (\ref{cheNS}) becomes as follows:
\begin{equation} \label{cheNS7}
     \lambda^{27} + 35 \lambda^{19} + 3 \lambda^{18} - 28 \lambda^{17} + 105 \lambda^{10} - 21 \lambda^{8} + 4410 = 0
\end{equation}
The characteristic equation (\ref{cheNS7}) has unstable roots $\lambda_{1,2}=1.2928 \pm 0.2032j$,
$|\mbox{arg}(\lambda_{1,2})| = 0.1560$ and therefore system (\ref{ChenSystem})
satisfy the necessary condition for exhibiting a double scroll attractor. Numerical simulation
of the system (\ref{ChenSystem}) for initial conditions
$(-9, -5, 14)$ is depicted in Fig.~\ref{Chenattr}.
\begin{figure}[ht]
\centering
\noindent
  \includegraphics[width=4.4in]{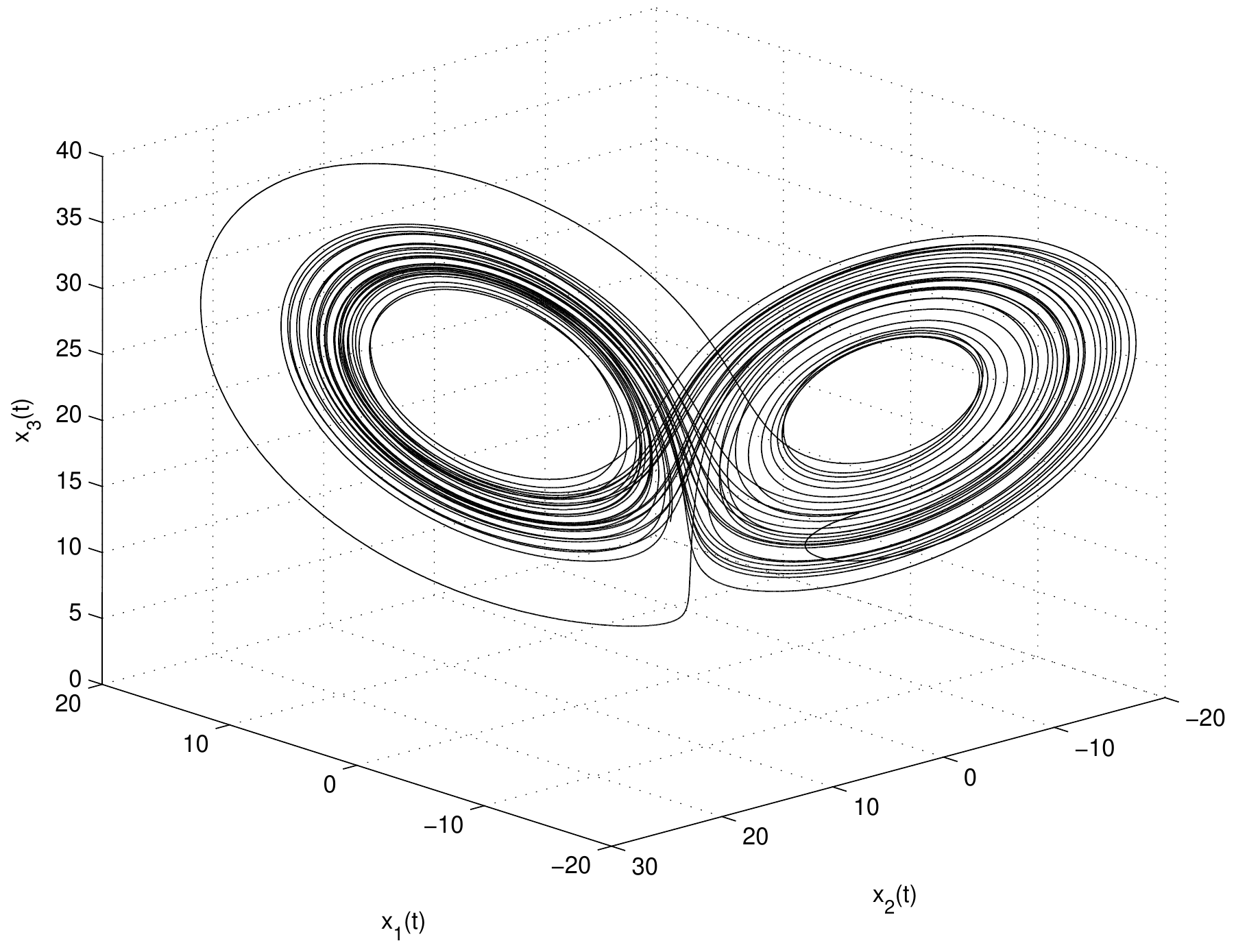}
   \caption{Double scroll attractor of Chen system (\ref{ChenSystem}) projected into 3D state space for $30\,sec$.}
   \label{Chenattr}
\end{figure}

\section{Conclusions}

In this paper we have presented the definitions for internal and external stability condition
of certain class of the linear and nonlinear fractional order system of finite dimension given in
state space, FODE or transfer function representation (polynomial).
It is important to note that stability and asymptotic behavior of fractional order system is not exponential type \cite{Bellman}
but it is in form of power law $t^{-\alpha}\, (\alpha \in R)$,  so called long memory behavior \cite{Matignon96}.

The results presented in this article are also applicable in robust stability investigation \cite{Chen2, Petras3, Petras4, Petras5},
stability of delayed system \cite{Deng, Chen1} and stability of discrete fractional order system \cite{Sierociuk, Matignon}.
Investigation of the fractional incommensurate order systems in state space, where space is deformed
by various order of derivatives in various directions is still open.

\section*{Acknowledgment}
This work was supported in part by the Slovak Grant Agency for Science
under grants VEGA: 1/3132/06, 1/4058/07, 1/0404/08, and APVV-0040-07.

\end{document}